\documentclass{article}

\newcommand{\email}[1]{}

\newenvironment{keywords}{\textbf{keywords}:}{\vspace*{1mm}\newline}
\newenvironment{AMS}{\textbf{AMS}:}{\vspace*{1mm}\newline}

\usepackage[margin=2cm]{geometry}
\usepackage [english]{babel}
\usepackage [english=british]{csquotes} % or american

\newif\ifportable
\portablefalse
\portabletrue

\usepackage{amssymb}

\usepackage{capt-of} % for \captionof

\usepackage[T1]{fontenc}
\usepackage{lmodern}
\usepackage{textcomp}
  
\usepackage{multirow}

\usepackage{algorithm}
\usepackage[noend]{algpseudocode}

\usepackage{booktabs}

\usepackage[	
		style=numeric-comp,
		hyperref=true,
		url=false,
		doi=false,
		isbn=false,
		eprint=false,
		firstinits=true,
		uniquename=init,
		maxbibnames=10,
		backend=bibtex,
		mincitenames=1,
		maxcitenames=3,
		]{biblatex}

\AtEveryBibitem{\clearfield{month}}
\AtEveryCitekey{\clearfield{month}}
\renewbibmacro{in:}{}

  \newtheorem{theorem}{Theorem}

  \newtheorem{definition}[theorem]{Definition}
  \newtheorem{remark}[theorem]{Remark}

\newif\ifcolouredsymbols

\colouredsymbolsfalse

\ifcolouredsymbols
  \newcommand{\mathcolour}{\color{blue}}

  \newcommand{\newmathcommand}[3][0]{
    \newcommand{#2}[#1]{\ensuremath{{\mathcolour#3}}}}
  \newcommand{\renewmathcommand}[3][0]{
    \renewcommand{#2}[#1]{\ensuremath{{\mathcolour#3}}}}
\else
  \newcommand{\mathcolour}{}

  \newcommand{\newmathcommand}[3][0]{
    \newcommand{#2}[#1]{\ensuremath{#3}}}
  \newcommand{\renewmathcommand}[3][0]{
    \renewcommand{#2}[#1]{\ensuremath{#3}}}
\fi

\newcommand{\MRIimage}{v}
\newcommand{\Reg}{\mathcal R}

\newmathcommand{\ParamReg}{\alpha}

\newmathcommand{\ParamThresh}{\eta}

\usepackage{mathtools}
\newmathcommand{\De}{\coloneqq}

\newmathcommand{\vecA}{x}
\newcommand{\vecB}{y}

\newmathcommand{\Grad}{\nabla}
\newmathcommand{\Div}{\operatorname{div}}
	
\newmathcommand{\absL}{|}
\newmathcommand{\absR}{|}
\newcommand{\abs}[1]{\absL#1\absR}

\newmathcommand{\normL}{\|}
\newmathcommand{\normR}{\|}
\newcommand{\norm}[1]{\normL#1\normR}

\renewmathcommand[2]{\angle}{\langle #1, #2 \rangle}
\DeclareMathOperator{\ArgminSymb}{{\mathcolour\operatorname{argmin}}}
\newcommand{\Argmin}[2]{\underset{#1}{\ArgminSymb} \;#2}

\DeclareMathOperator{\MinimizerSymb}{{\mathcolour\sharp}}
\newcommand{\Minimizer}[1]{#1^{\MinimizerSymb}}

\newmathcommand{\iterA}{k}
\newmathcommand{\IterA}{K}

\newmathcommand{\SubGradA}{p}
\newmathcommand{\SubGradB}{q}

\DeclareMathOperator{\ProxSymb}{{\mathcolour\operatorname{prox}}}
\newcommand{\Prox}[2]{\ProxSymb_{#1}(#2)}

\newcommand{\FunFromTo}[3]{#1\colon#2\rightarrow#3}

\newmathcommand{\Lin}{\operatorname{Lin}}

\newmathcommand{\OpA}{A}

\DeclareMathOperator{\AdjSymb}{{\mathcolour\ast}}
\newcommand{\Adj}[1]{#1^{\AdjSymb}}

\newmathcommand{\ProjA}{\mathcal P}
\newmathcommand{\id}{id}

\newmathcommand{\Id}{\mathcal I}

\newmathcommand{\N}{\mathbb N}
\newmathcommand{\R}{\mathbb R}
\newmathcommand{\RI}{\R_\infty}
\renewmathcommand{\RN}{[0,\infty)}
\newmathcommand{\C}{\mathbb C}

\renewmathcommand{\Reg}{J}

\newmathcommand{\MatAnisotropy}{\mathcal D}

\newmathcommand{\FGPstepsize}{s}
\newmathcommand{\FGPproxpoint}{y}

\newmathcommand{\ADMMrho}{\rho}
\newmathcommand{\ADMMlagrange}{L}
\newmathcommand{\ADMMlagmultX}{\mu}
\newmathcommand{\ADMMlagmultU}{\nu}
\newmathcommand{\ADMMvarX}{x}
\newmathcommand{\ADMMvarZ}{z}

\newmathcommand{\imaginary}{\mathrm{i}}
\newmathcommand{\MRIop}{\mathcal A}
\newmathcommand{\FT}{\mathcal F}
\newmathcommand{\Sampling}{\mathcal S}
\renewmathcommand{\Re}{\operatorname{Re}}
\newmathcommand{\AdjRe}{\Adj{\Re}}

\renewmathcommand{\MRIimage}{u}
\newmathcommand{\MRIsinfo}{v}
\newmathcommand{\MRIdata}{b}

\newmathcommand{\MRIsamplingSequence}{\pi}

\newmathcommand{\MRIdimImage}{N}
\newmathcommand{\MRIdimData}{M}

\newmathcommand{\InImage}{n}
\newmathcommand{\InData}{m}

\newmathcommand{\SymbolTV}{\operatorname{TV}}
\newmathcommand{\SymbolWTV}{\operatorname{wTV}}
\newmathcommand{\SymbolAQPL}{\operatorname{dTV}}

\newmathcommand{\WTVweight}{w}
\newmathcommand{\AQPLdirection}{\xi}

\newmathcommand{\GradSpace}{\mathbb G}
\newmathcommand{\MatSpace}{\mathbb M}

\newmathcommand{\UnitBall}{\mathbb U}
\newmathcommand{\ContraintSetImage}{\mathbb T}

\newcommand{\ie}{i.e.}
\newcommand{\eg}{e.g.}
\newcommand{\cf}{cf.}
\newcommand{\st}{s.t.}

\newcommand{\Quote}[1]{\enquote{#1}}

\newcommand{\TextData}[1]{\texttt{#1}}
\newcommand{\NameSheppLogan}{\TextData{SheppLogan}}
\newcommand{\NameBrainWeb}{\TextData{BrainWeb}}
\newcommand{\NameBrainWebA}{\TextData{BrainWebA}}
\newcommand{\NameBrainWebB}{\TextData{BrainWebB}}
\newcommand{\NameBrainWebC}{\TextData{BrainWebC}}
\newcommand{\NameNinon}{\TextData{patient}}
\newcommand{\NameNinonA}{\TextData{patientA}}
\newcommand{\NameNinonB}{\TextData{patientB}}

\newcommand{\FigureSettings}{
\centering
\footnotesize
}

\newcommand{\Portable}[2]{
\ifportable
\includegraphics{#1}
\vspace*{-3mm}
\else
\tikzsetnextfilename{#1}%
#2
\fi
}
\usepackage{tikz}
\newif\ifexternalize
\externalizetrue
\ifexternalize
  \usetikzlibrary{external}
 	\tikzset{external/system call= {pdflatex \tikzexternalcheckshellescape -halt-on-error -interaction=batchmode -jobname "\image" "\texsource"}} 
\fi

\newcommand{\Folder}{.}
\newlength{\Width}%
\newlength{\Height}%
\newlength{\Vspace}%
\newlength{\Hspace}%
\newlength{\x}%
\newlength{\y}%

\def\BarWidth{\Width}%
\newlength{\BarHeight}%
\newcommand{\NextPic}{%
  \addtolength{\x}{\Width}%
  \addtolength{\x}{\Hspace}}%
\newcommand{\NextRow}{%
  \setlength{\x}{0cm}%
  \addtolength{\y}{-\Height}%
  \addtolength{\y}{-\Vspace}}%
\newcommand{\FirstPic}{%
  \setlength{\x}{0cm}%
  \setlength{\y}{0cm}}%
\newcommand{\DrawText}[3]{\draw[color=black] (#1,#2) node [anchor=center] {#3};}%
\newcommand{\DrawVText}[3]{\draw[color=black] (#1,#2) node [anchor=center] {\rotatebox{90}{#3}};}%
\newcommand{\DrawColourBar}[8]{%
  \draw (#1,#2) node [anchor=north west] {\includegraphics[width=#8,height=#7,angle=#3]{\Folder/#6}};%
	\draw (#1+1.5mm,#2-.7*#8) node [anchor=west] {#4};% lower border
	\draw (#1+#7+1mm,#2-.7*#8) node [anchor=east] {#5};}% upper border
\begin{document}
% 
% \address{Centre for Medical Image Computing, University College London, London WC1E 6BT, UK}
%
% \ead{}
%
% \ams{xxx} % http://www.ams.org/msc/msc2010.html
% \pacs{xxx} % http://www.aip.org/pacs
% \submitto{\IP}
% 
% 
%
\title{Multi-Contrast MRI Reconstruction with Structure-Guided Total Variation\thanks{This research was funded by the EPSRC grants EP/H0464110/1 and EP/K009745/1 and the UCL Department of Computer Science.}}
\author{Matthias~J.~Ehrhardt\footnotemark[2] \and Marta~M.~Betcke\footnotemark[2]}
\maketitle
\renewcommand{\thefootnote}{\fnsymbol{footnote}}
\footnotetext[2]{Centre for Medical Image Computing, University College London, London WC1E 6BT, UK. (\url{matthias.ehrhardt.11@ucl.ac.uk}, \url{m.betcke@ucl.ac.uk}).}
% \footnotetext[3]{}
% \footnotetext[4]{}
%
\renewcommand{\thefootnote}{\arabic{footnote}}
\begin{abstract}
Magnetic resonance imaging (MRI) is a versatile imaging technique that allows different contrasts depending on the acquisition parameters. Many clinical imaging studies acquire MRI data for more than one of these contrasts---such as for instance $T_1$ and $T_2$ weighted images---which makes the overall scanning procedure very time consuming. As all of these images show the same underlying anatomy one can try to omit unnecessary measurements by taking the similarity into account during reconstruction. We will discuss two modifications of total variation---based on i) location and ii) direction---that take structural a priori knowledge into account and reduce to total variation in the degenerate case when no structural knowledge is available. We solve the resulting convex minimization problem with the alternating direction method of multipliers that separates the forward operator from the prior. For both priors the corresponding proximal operator can be implemented as an extension of the fast gradient projection method on the dual problem for total variation. We tested the priors on six data sets that are based on phantoms and real MRI images. In all test cases exploiting the structural information from the other contrast yields better results than separate reconstruction with total variation in terms of standard metrics like peak signal-to-noise ratio and structural similarity index. Furthermore, we found that exploiting the two dimensional directional information results in images with well defined edges, superior to those reconstructed solely using a priori information about the edge location.
\end{abstract}
%\begin{AMS} \added[id=ME]{47A52, Ill-posed problems, regularization}, \added[id=ME]{65J22, Inverse problems}, \added[id=ME]{49M30, Other methods (Calculus of variations and optimal control; optimization)} , \added[id=ME]{94A08, Image processing (compression, reconstruction, etc.) }  \end{AMS}
% 
\pagestyle{myheadings}
\thispagestyle{plain}
\markboth{M. J. EHRHARDT, M. M. BETCKE}{MULTI-CONTRAST MRI WITH STRUCTURE-GUIDED TV}
%
%STRUCTURE-GUIDED MULTI-CONTRAST MRI 
%MULTI-CONTRAST MRI WITH STRUCTURE-GUIDED TV
%MULTI-CONTRAST MRI WITH STRUCTURE-GUIDED TOTAL VARIATION
%MULTI-CONTRAST MRI RECONSTRUCTION WITH STRUCTURE-505050GUIDED TOTAL VARIATION
%10101010102020202020303030303040404040405050505050

\section{Introduction}
\subsection{Multi-Contrast Magnetic Resonance Imaging}
Magnetic resonance imaging (MRI) is a well established imaging modality with numerous applications. One of its key advantages is versatility: depending on image acquisition protocol, images with very different contrast and informational content can be acquired \cite{Liang2000,McRobbie2006}. Most common are images that are weighted by the relaxation times $T_1$ and $T_2$ but many more options are available. In clinical applications, often not one but several MRI images with different contrasts are acquired during one session. As an example, the UK Biobank\footnote{\url{http://biobank.ctsu.ox.ac.uk/crystal/label.cgi}, accessed August 14, 2015} contains for each subject MRI data for images weighted not only by $T_1$ and $T_2$ but also for images that are fluid-suppressed (FLAIR), or show susceptibility, diffusion or function. All of these data have to be acquired sequentially one at a time, which makes the whole scanning procedure rather lengthy. 
Therefore, shortening the acquisition time would not only reduce patient's discomfort but would increase the patient throughput leading to more efficient use of the scanning facilities.
\subsection{Magnetic Resonance Imaging and Compressed Sensing}
\newcommand{\ColourAQPL}{Paired-6-6}
\newcommand{\MarkAQPL}{square*}
\newcommand{\MarkAngleAQPL}{0}
\newcommand{\MarkSizeAQPL}{.8pt}

\newcommand{\ColourTV}{Paired-6-4}
\newcommand{\MarkTV}{*}
\newcommand{\MarkAngleTV}{0}
\newcommand{\MarkSizeTV}{1pt}

\newcommand{\ColourES}{RdBu-11-3}
\newcommand{\MarkES}{x}
\newcommand{\MarkAngleES}{0}
\newcommand{\MarkSizeES}{1pt}

\newcommand{\ColourWTV}{Paired-6-2}
\newcommand{\MarkWTV}{diamond*}
\newcommand{\MarkAngleWTV}{0}
\newcommand{\MarkSizeWTV}{1.2pt}

\newlength{\WidthSmall}%
\newlength{\WidthLarge}%

\renewcommand{\DrawColourBar}[8]{%
  \draw (#1,#2) node [anchor=south west] {\includegraphics[width=#8,height=#7,angle=#3,frame]{\Folder/#6}};%
	\draw (#1+#8+1mm,#2+1mm) node [anchor=south west] {#4};% lower border
	\draw (#1+#8+1mm,#2+#7+2mm) node [anchor=north west] {#5};}% upper border

\newcommand{\DrawColourBarPhaseRotated}[5]{%
		\begin{scope}%
		\clip (#1,#2) circle (#3);%
		\draw(#1,#2) node {\includegraphics[width=#4,rotate=90]{#5}};%
		\draw(#1,#2) circle (#3);
		\end{scope}
		\draw(#1,#2+#3) node[anchor=south] {$0$};%
		\draw(#1,#2-#3) node[anchor=north] {$\pi$};}%

\newcommand{\DrawColourBarPhase}[5]{%
		\DrawColourBarPhaseWO{#1}{#2}{#3}{#4}{#5}
		\draw(#1+#3,#2) node[anchor=west] {$0$};%
		\draw(#1-#3,#2) node[anchor=east] {$\pi$};}%
		
\newcommand{\DrawColourBarPhaseWO}[5]{%
		\begin{scope}%
		\clip (#1,#2) circle (#3);%
		\draw(#1,#2) node {\includegraphics[width=#4]{#5}};%
		\draw(#1,#2) circle (#3);
		\end{scope}}%

\newcommand{\FolderBase}{pics}

\setlength{\Width}{2.0cm}%
\setlength{\Height}{2.0cm}%
\setlength{\Vspace}{.0cm}%
\def\BarWidth{\Width}%
\def\BarHeight{3mm}%

\newcommand{\DrawIm}[1]{%
  \draw[] (\x,\y) node[anchor=south west] {\includegraphics[height=\Height]{#1}};}%
  
\newcommand{\DataSetA}{phantom}
\newcommand{\DataSetDateA}{20150729_1551}
\newcommand{\FolderDataA}{\FolderBase/data/MRI2_Jun15_\DataSetA_\DataSetDateA}
\newcommand{\DrawImA}[1]{\draw[] (\x,\y) node[anchor=south west] {\includegraphics[height=\Height]{#1}};}%
\newcommand{\DrawZoomA}[1]{\draw[] (\x,\y) node[anchor=south west] {\includegraphics[clip,trim=4cm 5cm 2cm 1cm,height=\Height]{#1}};}%

\newcommand{\DataSetB}{BrainWebNormalTransverse1}
\newcommand{\DataSetDateB}{20150724_1642}
\newcommand{\FolderDataB}{\FolderBase/data/MRI2_Jun15_\DataSetB_\DataSetDateB}
\newcommand{\DrawImB}[1]{\draw[] (\x,\y) node[anchor=south west] {\includegraphics[width=\Height,clip,trim=.5cm .5cm .5cm .5cm,angle=90]{#1}};}%
\newcommand{\DrawZoomB}[1]{\draw[] (\x,\y) node[anchor=south west] {\includegraphics[clip,trim=5cm 2cm 1cm 4cm,width=\Height,angle=90]{#1}};}%

\newcommand{\DataSetC}{BrainWebNormalTransverse2}
\newcommand{\DataSetDateC}{20150724_1638}
\newcommand{\FolderDataC}{\FolderBase/data/MRI2_Jun15_\DataSetC_\DataSetDateC}
\newcommand{\DrawImC}[1]{\draw[] (\x,\y) node[anchor=south west] {\includegraphics[width=\Height,clip,trim=.5cm .5cm .5cm .5cm,angle=90]{#1}};}%
\newcommand{\DrawZoomC}[1]{\draw[] (\x,\y) node[anchor=south west] {\includegraphics[width=\Height,clip,trim=1.5cm 3cm 5cm 3.5cm,angle=90]{#1}};}%

\newcommand{\DataSetD}{BrainWebNormalCoronal1}
\newcommand{\DataSetDateD}{20150724_1644}
\newcommand{\FolderDataD}{\FolderBase/data/MRI2_Jun15_\DataSetD_\DataSetDateD}
\newcommand{\DrawImD}[1]{\draw[] (\x,\y) node[anchor=south west] {\includegraphics[height=\Height,clip,trim=1.1cm 0cm 1.1cm 2.2cm,angle=0]{#1}};}%
\newcommand{\DrawZoomD}[1]{\draw[] (\x,\y) node[anchor=south west] {\includegraphics[width=\Height,clip,trim=1.5cm 3cm 5cm 3.5cm,angle=0]{#1}};}%

\newcommand{\DataSetE}{NinonTransverse}
\newcommand{\DataSetDateE}{20150727_1739}
\newcommand{\FolderDataE}{\FolderBase/data/MRI2_Jun15_\DataSetE_\DataSetDateE}
\newcommand{\DrawImE}[1]{\draw[] (\x,\y) node[anchor=south west] {\includegraphics[height=\Height,clip,trim=.5cm .5cm .5cm .5cm,angle=0]{#1}};}%
\newcommand{\DrawZoomEA}[1]{\draw[] (\x,\y) node[anchor=south west] {\includegraphics[width=\Height,clip,trim=4cm 1cm 1cm 4cm,angle=0]{#1}};}%
\newcommand{\DrawZoomEB}[1]{\draw[] (\x,\y) node[anchor=south west] {\includegraphics[width=\Height,clip,trim=1.5cm 3.5cm 4cm 2cm,angle=0]{#1}};}%
\newcommand{\DrawZoomEC}[1]{\draw[] (\x,\y) node[anchor=south west] {\includegraphics[width=\Height,clip,trim=4cm 1cm 1cm 4cm,angle=0]{#1}};}%

\newcommand{\DataSetF}{NinonSagital}
\newcommand{\DataSetDateF}{20150727_1742}
\newcommand{\FolderDataF}{\FolderBase/data/MRI2_Jun15_\DataSetF_\DataSetDateF}
\newcommand{\DrawImF}[1]{\draw[] (\x,\y) node[anchor=south west] {\includegraphics[height=\Height,clip,trim=0cm 0cm 0cm 0cm,angle=0]{#1}};}%
\newcommand{\DrawZoomF}[1]{\draw[] (\x,\y) node[anchor=south west] {\includegraphics[width=\Height,clip,trim=1cm 3cm 4cm 2cm,angle=0]{#1}};}%

\newcommand{\DrawGT}[1]{
	\addtolength{\x}{2mm}
  \DrawImLocal{\FolderBase/data/MRI2_Jun15_\DataSetName_groundtruth_#1}}

\newcommand{\DrawGTZoom}[1]{
	\addtolength{\x}{2mm}
  \DrawZoomLocal{\FolderBase/data/MRI2_Jun15_\DataSetName_groundtruth_#1}}

\newcommand{\DrawBarSSIM}{
\DrawColourBar{-\Width+7mm}{\y+1mm}{180}{$\leq 0$}{1}{pics/results/colormap_ssim}{\BarWidth-4mm}{\BarHeight}}

\newcommand{\DrawBarError}{
\DrawColourBar{-\Width+7mm}{\y+1mm}{180}{$\leq -0.3$}{$\geq +0.3$}{pics/results/colormap_error}{\BarWidth-4mm}{\BarHeight}}

\newcommand{\DrawRow}[3]{
  \DrawImLocal{\FolderBase/results/#1_early_stopping_\DataSetName_\SamplingName_best_#3}%
  \NextPic \DrawImLocal{\FolderBase/results/#1_TV-FGP_\DataSetName_\SamplingName_best_#3}%
  \NextPic \DrawImLocal{\FolderBase/results/#1_wTV-FGP_#2gt_\DataSetName_\SamplingName_best_#3}%
  \NextPic \DrawImLocal{\FolderBase/results/#1_AQPL-FGP_#2gt_\DataSetName_\SamplingName_best_#3}}

\newcommand{\DrawRowZoom}[3]{
  \DrawZoomLocal{\FolderBase/results/#1_early_stopping_\DataSetName_\SamplingName_best_#3}%
  \NextPic \DrawZoomLocal{\FolderBase/results/#1_TV-FGP_\DataSetName_\SamplingName_best_#3}%
  \NextPic \DrawZoomLocal{\FolderBase/results/#1_wTV-FGP_#2gt_\DataSetName_\SamplingName_best_#3}%
  \NextPic \DrawZoomLocal{\FolderBase/results/#1_AQPL-FGP_#2gt_\DataSetName_\SamplingName_best_#3}}

\newcommand{\DrawSampling}{
    \FirstPic    
    \setlength{\Width}{\WidthSmall}%
    \setlength{\Height}{\WidthSmall}%
    \addtolength{\x}{-\WidthSmall}
    \addtolength{\x}{-1mm}
    \addtolength{\y}{-\WidthSmall}
    \addtolength{\y}{\WidthLarge}
    \DrawIm{\FolderDataLocal/MRI2_Jun15_\DataSetName_sampling_\SamplingName}%
    \DrawText{\x+.5*\WidthSmall+1mm}{\y-1mm}{sampling}}

\newcommand{\DrawRowLabels}{
    \setlength{\x}{0.5\Width + 1mm}
		\addtolength{\y}{-1mm}
    \DrawText{\x}{\y}{no prior}
    \NextPic\DrawText{\x}{\y}{\phantom{g}$ \SymbolTV $\phantom{f}}
    \NextPic\DrawText{\x}{\y}{\phantom{g}$ \SymbolWTV $\phantom{f}}
    \NextPic\DrawText{\x}{\y}{\phantom{g}$ \SymbolAQPL $\phantom{f}}
		\addtolength{\x}{2mm}
		\NextPic\DrawText{\x}{\y}{ground truth}}

\newcommand{\SetLengths}{
  \setlength{\WidthSmall}{2.0cm}%
  \setlength{\WidthLarge}{2.2cm}%
  \setlength{\Vspace}{.0cm}%
  \setlength{\Hspace}{-.1cm}%
  \setlength{\Width}{\WidthLarge}%
  \setlength{\Height}{\WidthLarge}}

\newcommand{\DrawResults}{
  \begin{tikzpicture}%
    \SetLengths
    \FirstPic\DrawRow{T1}{T2}{solution} 
		\NextPic \DrawGT{T1} 
		\NextRow \DrawRowZoom{T1}{T2}{solution} 
		\NextPic \DrawGTZoom{T1} 
    \addtolength{\y}{-1mm}
    \NextRow \DrawRow{T2}{T1}{solution}
		\NextPic \DrawGT{T2} 
    \NextRow \DrawRowZoom{T2}{T1}{solution}
		\NextPic \DrawGTZoom{T2} 
    \DrawRowLabels
    \DrawSampling
  \end{tikzpicture}}
  
\newcommand{\DrawResultsTone}{
  \begin{tikzpicture}%
    \SetLengths
    \FirstPic \DrawRow{T1}{T2}{solution}
    \NextPic \DrawGT{T1} 
		\NextRow \DrawRowZoom{T1}{T2}{solution}
		\NextPic \DrawGTZoom{T1} 
    \DrawRowLabels
    \DrawSampling
  \end{tikzpicture}}
  
\newcommand{\DrawResultsTtwo}{
  \begin{tikzpicture}%
    \SetLengths
    \FirstPic \DrawRow{T2}{T1}{solution}
		\NextPic \DrawGT{T2} 
    \NextRow \DrawRowZoom{T2}{T1}{solution}
		\NextPic \DrawGTZoom{T2} 
    \DrawRowLabels
    \DrawSampling
  \end{tikzpicture}}
  
\newcommand{\DrawResultsSSIM}{
  \begin{tikzpicture}%
    \SetLengths
    \FirstPic \DrawRow{T1}{T2}{solution}
		\NextPic \DrawGT{T1} 
    \NextRow \DrawRow{T1}{T2}{ssimmap}
    \DrawBarSSIM
    \addtolength{\y}{-1mm}
    \NextRow \DrawRow{T2}{T1}{solution}
		\NextPic \DrawGT{T2} 
    \NextRow \DrawRow{T2}{T1}{ssimmap}
    \DrawRowLabels
    \DrawSampling
  \end{tikzpicture}}
  
\newcommand{\DrawResultsToneSSIM}{
  \begin{tikzpicture}%
    \SetLengths
    \FirstPic \DrawRow{T1}{T2}{solution}
		\NextPic \DrawGT{T1} 
    \NextRow \DrawRow{T1}{T2}{ssimmap}
    \DrawBarSSIM
    \DrawRowLabels
    \DrawSampling
  \end{tikzpicture}}
  
\newcommand{\DrawResultsError}{
  \begin{tikzpicture}%
    \SetLengths
    \FirstPic \DrawRow{T1}{T2}{solution}
		\NextPic \DrawGT{T1} 
    \NextRow \DrawRow{T1}{T2}{error}
    \addtolength{\y}{-1mm}
    \NextRow \DrawRow{T2}{T1}{solution}
		\NextPic \DrawGT{T2} 
    \NextRow \DrawRow{T2}{T1}{error}
    \DrawBarError
    \DrawRowLabels
    \DrawSampling
  \end{tikzpicture}}
  
\newcommand{\DrawResultsToneError}{
  \begin{tikzpicture}%
    \SetLengths
    \FirstPic \DrawRow{T1}{T2}{solution}
		\NextPic \DrawGT{T1} 
    \NextRow \DrawRow{T1}{T2}{error}
    \DrawBarError
    \DrawRowLabels
    \DrawSampling
  \end{tikzpicture}}
  
\newcommand{\DrawResultsTtwoError}{
  \begin{tikzpicture}%
    \SetLengths
    \FirstPic \DrawRow{T2}{T1}{solution}
		\NextPic \DrawGT{T2} 
    \NextRow \DrawRow{T2}{T1}{error}
    \DrawBarError
    \DrawRowLabels
    \DrawSampling
  \end{tikzpicture}}

% 	\NextPic \DrawImLocal{\FolderBase/results/T1_early_stopping_\DataSetName_\SamplingName_best_error}%
% 	\DrawVText{\x - 1mm}{\y + .5\Width + 1mm}{error T1}
% 	\NextPic \DrawImLocal{\FolderBase/results/T1_TV-FGP_\DataSetName_\SamplingName_best_error}%
% 	\NextPic \DrawImLocal{\FolderBase/results/T1_wTV-FGP_T2gt_\DataSetName_\SamplingName_best_error}%
% 	\NextPic \DrawImLocal{\FolderBase/results/T1_AQPL-FGP_T2gt_\DataSetName_\SamplingName_best_error}%}{\BarHeight}
% 	\NextRow%
% 	\addtolength{\x}{5mm}
% 	\NextPic \DrawImLocal{\FolderBase/results/T2_early_stopping_\DataSetName_\SamplingName_best_solution}%
% 	\DrawVText{\x - 1mm}{\y + .5\Width + 1mm}{T2}
% 	\NextPic \DrawImLocal{\FolderBase/results/T2_TV-FGP_\DataSetName_\SamplingName_best_solution}%
% 	\NextPic \DrawImLocal{\FolderBase/results/T2_wTV-FGP_T1gt_\DataSetName_\SamplingName_best_solution}%
% 	\NextPic \DrawImLocal{\FolderBase/results/T2_AQPL-FGP_T1gt_\DataSetName_\SamplingName_best_solution}%
% 	\NextRow%
% 	\addtolength{\x}{5mm}
% 	\NextPic \DrawImLocal{\FolderBase/results/T2_early_stopping_\DataSetName_\SamplingName_best_error}%
% 	\DrawVText{\x - 1mm}{\y + .5\Width + 1mm}{error T2}
% 	\NextPic \DrawImLocal{\FolderBase/results/T2_TV-FGP_\DataSetName_\SamplingName_best_error}%
% 	\NextPic \DrawImLocal{\FolderBase/results/T2_wTV-FGP_T1gt_\DataSetName_\SamplingName_best_error}%
% 	\NextPic \DrawImLocal{\FolderBase/results/T2_AQPL-FGP_T1gt_\DataSetName_\SamplingName_best_error}%	
% \end{tikzpicture}}%

\begin{figure}%
\FigureSettings
\Portable{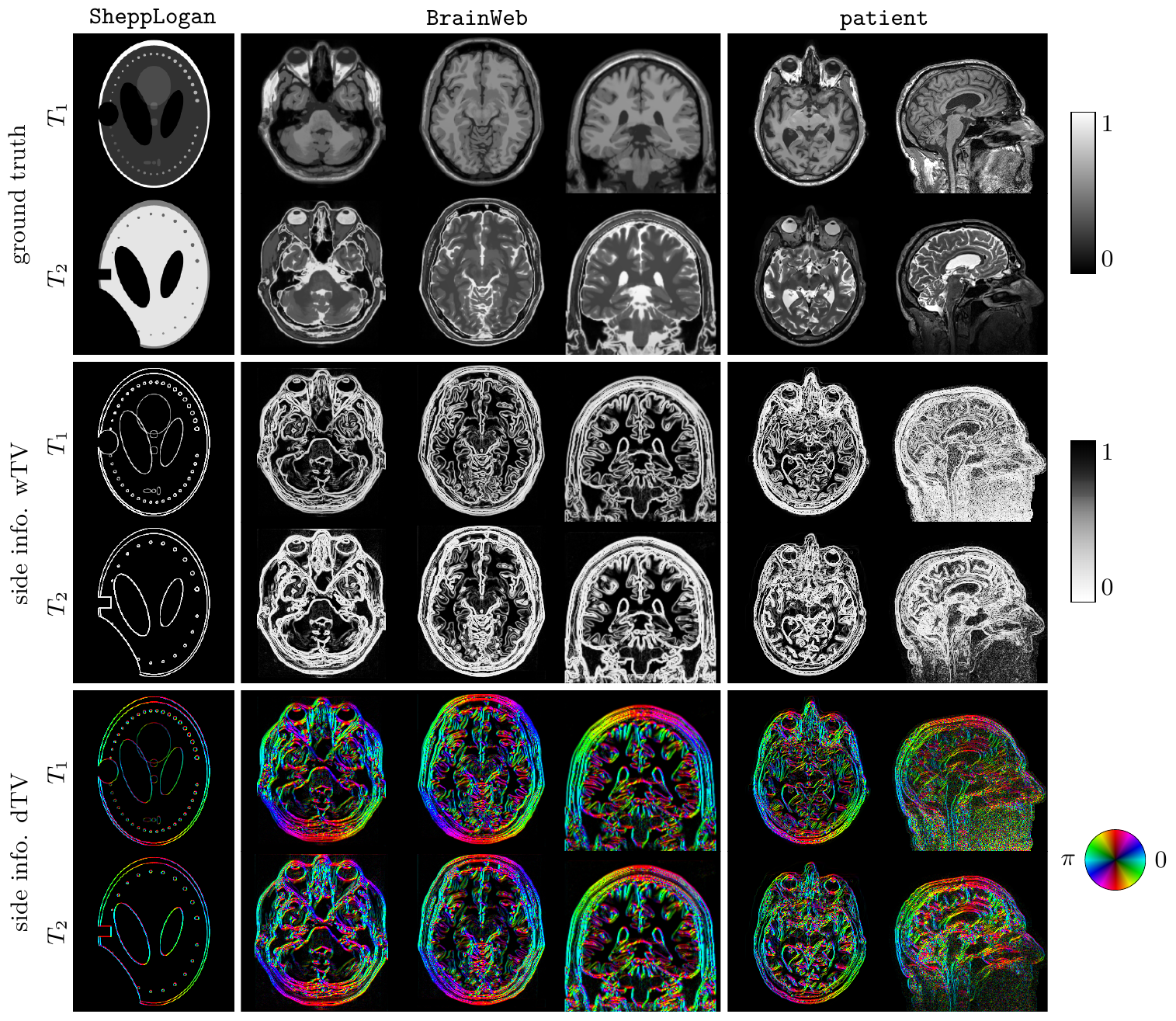}{
\setlength{\Width}{2.1cm}%
\setlength{\Height}{2.1cm}%
\setlength{\Vspace}{.0cm}%
\setlength{\Hspace}{-.1mm}%
\begin{tikzpicture}%
\FirstPic\DrawImA{\FolderDataA/MRI2_Jun15_\DataSetA_groundtruth_T1}%
\DrawVText{\x - 6mm}{\y + 1mm}{ground truth}
\DrawVText{\x - 1mm}{\y + .5\Width + 1mm}{$T_1$}
\DrawText{\x + 0.5\Width + 1mm}{\y + \Width + 3mm}{\NameSheppLogan}
\addtolength{\x}{1mm}
\NextPic \DrawImB{\FolderDataB/MRI2_Jun15_\DataSetB_groundtruth_T1}%
\NextPic \DrawImC{\FolderDataC/MRI2_Jun15_\DataSetC_groundtruth_T1}%
\DrawText{\x + 0.5\Width + 1mm}{\y + \Width + 3mm}{\phantom{f}\NameBrainWeb\phantom{g}}
\NextPic \DrawImD{\FolderDataD/MRI2_Jun15_\DataSetD_groundtruth_T1}%
\addtolength{\x}{1mm}
\NextPic \DrawImE{\FolderDataE/MRI2_Jun15_\DataSetE_groundtruth_T1}%
\DrawText{\x + \Width + 1mm}{\y + \Width + 3mm}{\phantom{f}\NameNinon\phantom{g}}
\NextPic \DrawImF{\FolderDataF/MRI2_Jun15_\DataSetF_groundtruth_T1}%
\NextRow \DrawImA{\FolderDataA/MRI2_Jun15_\DataSetA_groundtruth_T2}%
\DrawVText{\x - 1mm}{\y + .5\Width + 1mm}{$T_2$}
\addtolength{\x}{1mm}
\NextPic \DrawImB{\FolderDataB/MRI2_Jun15_\DataSetB_groundtruth_T2}%
\NextPic \DrawImC{\FolderDataC/MRI2_Jun15_\DataSetC_groundtruth_T2}%
\NextPic \DrawImD{\FolderDataD/MRI2_Jun15_\DataSetD_groundtruth_T2}%
\addtolength{\x}{1mm}
\NextPic \DrawImE{\FolderDataE/MRI2_Jun15_\DataSetE_groundtruth_T2}%
\NextPic \DrawImF{\FolderDataF/MRI2_Jun15_\DataSetF_groundtruth_T2}%
\addtolength{\x}{-4mm}
\DrawColourBar{\x+\Width+7mm}{\y+.5*\Height}{180}{0}{1}{pics/results/gray}{\BarWidth}{\BarHeight}
\addtolength{\y}{-1mm}
\NextRow \DrawImA{\FolderBase/results/\DataSetA_sideinfo_wTV_T1}%
\DrawVText{\x - 6mm}{\y + 1mm}{side info. $ \SymbolWTV $}
\DrawVText{\x - 1mm}{\y + .5\Width + 1mm}{$T_1$}
\addtolength{\x}{1mm}
\NextPic \DrawImB{\FolderBase/results/\DataSetB_sideinfo_wTV_T1}%
\NextPic \DrawImC{\FolderBase/results/\DataSetC_sideinfo_wTV_T1}%
\NextPic \DrawImD{\FolderBase/results/\DataSetD_sideinfo_wTV_T1}%
\addtolength{\x}{1mm}
\NextPic \DrawImE{\FolderBase/results/\DataSetE_sideinfo_wTV_T1}%
\NextPic \DrawImF{\FolderBase/results/\DataSetF_sideinfo_wTV_T1}%
\NextRow \DrawImA{\FolderBase/results/\DataSetA_sideinfo_wTV_T2}%
\DrawVText{\x - 1mm}{\y + .5\Width + 1mm}{$T_2$}
\addtolength{\x}{1mm}
\NextPic \DrawImB{\FolderBase/results/\DataSetB_sideinfo_wTV_T2}%
\NextPic \DrawImC{\FolderBase/results/\DataSetC_sideinfo_wTV_T2}%
\NextPic \DrawImD{\FolderBase/results/\DataSetD_sideinfo_wTV_T2}%
\addtolength{\x}{1mm}
\NextPic \DrawImE{\FolderBase/results/\DataSetE_sideinfo_wTV_T2}%
\NextPic \DrawImF{\FolderBase/results/\DataSetF_sideinfo_wTV_T2}%
\addtolength{\x}{-4mm}
\DrawColourBar{\x+\Width+7mm}{\y+.5*\Height}{0}{0}{1}{pics/results/sideinfo_wTV_colormap}{\BarWidth}{\BarHeight}
\addtolength{\y}{-1mm}
\NextRow \DrawImA{\FolderBase/results/\DataSetA_sideinfo_AQPL_T1}%
\DrawVText{\x - 6mm}{\y + 1mm}{side info. $ \SymbolAQPL $}
\DrawVText{\x - 1mm}{\y + .5\Width + 1mm}{$T_1$}
\addtolength{\x}{1mm}
\NextPic \DrawImB{\FolderBase/results/\DataSetB_sideinfo_AQPL_T1}%
\NextPic \DrawImC{\FolderBase/results/\DataSetC_sideinfo_AQPL_T1}%
\NextPic \DrawImD{\FolderBase/results/\DataSetD_sideinfo_AQPL_T1}%
\addtolength{\x}{1mm}
\NextPic \DrawImE{\FolderBase/results/\DataSetE_sideinfo_AQPL_T1}%
\NextPic \DrawImF{\FolderBase/results/\DataSetF_sideinfo_AQPL_T1}%
\NextRow \DrawImA{\FolderBase/results/\DataSetA_sideinfo_AQPL_T2}%
\DrawVText{\x - 1mm}{\y + .5\Width + 1mm}{$T_2$}
\addtolength{\x}{1mm}
\NextPic \DrawImB{\FolderBase/results/\DataSetB_sideinfo_AQPL_T2}%
\NextPic \DrawImC{\FolderBase/results/\DataSetC_sideinfo_AQPL_T2}%
\NextPic \DrawImD{\FolderBase/results/\DataSetD_sideinfo_AQPL_T2}%
\addtolength{\x}{1mm}
\NextPic \DrawImE{\FolderBase/results/\DataSetE_sideinfo_AQPL_T2}%
\NextPic \DrawImF{\FolderBase/results/\DataSetF_sideinfo_AQPL_T2}%
\DrawColourBarPhase{\x+\Width+10mm}{\y+\Height + 0mm}{.4cm}{.8cm}{pics/results/hsv_magn}%
%
% \DrawColourBarPhaseWO{\x+\Width+10mm}{\y+0.2*\Height + 2mm}{.4cm}{2cm}{pics/results/hsv_pi_00}
% \DrawColourBarPhaseWO{\x+\Width+10mm}{\y+0.4*\Height + 2mm}{.4cm}{2cm}{pics/results/hsv_pi_02}
% \DrawColourBarPhaseWO{\x+\Width+10mm}{\y+0.6*\Height + 2mm}{.4cm}{2cm}{pics/results/hsv_pi_04}
% \DrawColourBarPhaseWO{\x+\Width+10mm}{\y+0.8*\Height + 2mm}{.4cm}{2cm}{pics/results/hsv_pi_05}
% \DrawColourBarPhaseWO{\x+\Width+10mm}{\y+1.0*\Height + 2mm}{.4cm}{2cm}{pics/results/hsv_pi_06}
% \DrawColourBarPhaseWO{\x+\Width+10mm}{\y+1.2*\Height + 2mm}{.4cm}{2cm}{pics/results/hsv_pi_08}
% \DrawColourBarPhaseWO{\x+\Width+10mm}{\y+1.4*\Height + 2mm}{.4cm}{2cm}{pics/results/hsv_pi_09}
% \DrawColourBarPhase{\x+\Width+10mm}{\y+1.6*\Height + 2mm}{.4cm}{2cm}{pics/results/hsv_pi_10}
% \draw[color=black,line width=.8mm,->] (\x+\Width+20mm,\y+3mm) -> (\x+\Width+20mm,\y+2*\Height-3mm);
% \draw node at (\x+\Width+23mm,\y+19mm) {\rotatebox{270}{\rotatebox{90}{1} \hspace*{3mm} magnitude \hspace*{3mm} \rotatebox{90}{0}}};
% \DrawText{\x+\Width+10mm}{\y+2*\Width}{phase}
\end{tikzpicture}%
}%
\caption{Ground truth $T_1$ and $T_2$ images with the side information that is exploited by weighted and directional total variation.} \label{FIG:GROUNDTRUTH}%
\end{figure}%
To speed up the scanning procedure, it has been proposed almost a decade ago to apply compressed sensing \cite{Candes2006a,Candes2006b,Donoho2006,Eldar2012} to MRI \cite{Lustig2007} which is still an active research topic \cite{Chang2006,Lingala2011,Knoll2011,Ravishankar2011,Bilgic2011,Chen2013,Ma2013,Huang2014,Wang2014,Tremoulheac2014,Sodickson2015}. One of the main ideas of compressed sensing is to acquire fewer measurements and to solve the reconstruction problem by exploiting a priori knowledge about the solution. Initially, the a priori knowledge has been sparsity in a wavelet basis and penalizing large total variations; the latter being related to sparsity of the image gradient. Over the years many other forms of a priori knowledge have been proposed for MRI reconstruction such as higher order total variation \cite{Knoll2011}, sparsity in a self-learned dictionary \cite{Ravishankar2011} or regularization of dynamic sequences with the nuclear norm \cite{Lingala2011,Tremoulheac2014} to name just a few. In a multi-contrast MRI scan, the images have very different information content, but as they are acquired from the same patient anatomy we know a priori that they are likely to show very similar structures \cite{Bilgic2011,Huang2014}. An example of a $T_1$ and $T_2$ weighted pair of MRI images of the same subject is shown on the right in figure \ref{FIG:GROUNDTRUTH}. Parallel MRI \cite{Pruessmann1999,Griswold2002a,Larkman2007,Uecker2014} is another example of image reconstruction problem which can benefit from exploiting common information. In \cite{Chen2013} joint reconstruction of different coil images is performed in the framework of compressed sensing.
\subsection{Contributions}
In this paper we aim to exploit the expected redundancy in a series of multi-contrast MRI images by extracting information about i) the location of edges and ii) the direction of edges from one contrast to aid the reconstruction of the other. We propose two priors that enable us to incorporate a priori structural knowledge into a total variation functional. In both cases the prior is convex such that we can use algorithms from convex optimization to solve the minimization problem. A double-split allows us to apply the alternating direction method of multipliers (ADMM) where all but one update are in closed-form. An extension of the fast gradient projection method first proposed for the standard total variation in \cite{Beck2009} is used to efficiently compute proximal operator for both priors.
\subsection{Related Work}
% \cite{Arridge2007}: Used smooth versions, only isotropic, hence similar to weighted TV
% \cite{Grasmair2009}: Adaptive regularization parameter selection - try to guess from the image, we use the other contrast in wTV; some anisotropy is mentioned, too
% \cite{Knoll2010}: denoising in MRI with adaptive regularization parameter
% \cite{Berkels2006}: do not have to cite this, rotation of the gradient and then l1 norm. not sure what to do.
% \cite{Bayram2012}: directional TV, denoising, rotate and scale
% \cite{Lenzen2015}, adaptive TV, solution driven, deblurring, inpainting
% \cite{Grasmair2010}, adaptive TV, structure tensor
% 
% \cite{Estellers2015}. adaptive TV, structure tensor
% \cite{Leahy1991,Bowsher2004,Vunckx2012,Arridge2007,Moeller2012,Fang2013,Kaipio1999}: one sided reconstruction in different applications, none in multi-contrast MRI; only \cite{Kaipio1999} directional, but not TV
% 
% \cite{Ehrhardt2014c,Ehrhardt2015}: similar to parallel level sets; used for joint reconstruction of PET-MRI and colour image processing.  We will point out in \S\ref{SEC:METHODS} the similarities and differences in more detail.
% 
In this work we propose extensions of total variation, based on i) location and ii) direction, that can exploit structural a-priori knowledge and apply it to the multi-contrast MRI setting where structural information is available from another contrast. In this context we group the related work into the four following classes:

\paragraph{Total Variation with Local Weighting}
Extensions of total variation or similar edge-preserving priors with spatially varying regularization parameter have been used before for optical tomography \cite{Arridge2007} and image denoising \cite{Grasmair2009,Knoll2010}. While the weights are a priori defined by side information in \cite{Arridge2007}, they are estimated based on local statistics in \cite{Grasmair2009,Knoll2010}. In that respect this contribution improves upon \cite{Arridge2007} as our algorithm can handle a non-smooth formulation of the prior.

\paragraph{Total Variation and Directional Information}
It has been proposed to include directional information into the total variation functional either by rotating the coordinate system and using locally the $\ell^1$-norm \cite{Berkels2006} or by scaling preferred directions and applying the $\ell^2$-norm \cite{Grasmair2010,Bayram2012,Estellers2015,Lenzen2015}. The directions are globally constant and predefined in \cite{Bayram2012} and based on the image content in \cite{Berkels2006,Grasmair2010,Estellers2015,Lenzen2015}. Our approach for directional information in the total variation functional shares similarities with \cite{Grasmair2010,Estellers2015,Lenzen2015}. While \cite{Grasmair2010,Estellers2015,Lenzen2015} compute the directions and scaling from the structure tensor of the current image estimate or the noisy input image, we project the gradient in the total variation functional onto a predefined vector field given by the other contrast. 

\paragraph{One-Sided Reconstruction}
Incorporating structural information by a prior has been used in other settings such as combined positron emission tomography (PET) with computed tomography or MRI \cite{Leahy1991,Bowsher2004,Vunckx2012}, optical tomography \cite{Arridge2007}, remote sensing \cite{Moeller2012,Fang2013} or electric impedance tomography \cite{Kaipio1999} but to the best of our knowledge has not been applied to multi-contrast MRI. In addition, only \cite{Arridge2007} and \cite{Kaipio1999} share similarities with our approach. In \cite{Arridge2007} the authors propose to locally adapt the weight of the prior isotropically and used a smoothed penalty function to facilitate diffusion techniques for reconstruction. On the other hand while the prior in \cite{Kaipio1999} is anisotropic, as it directionally dependent, it reduces to a quadratic prior when no edge information is available. In contrast, the here proposed priors reduce to total variation in absence of additional information.

\paragraph{Parallel Level Sets}
The directional extension of total variation is related to the idea of measuring the difference in structure of two images by means of parallel level sets. A symmetric version of the latter has been used for joint reconstruction of PET-MRI \cite{Ehrhardt2014,Ehrhardt2015} and colour image processing \cite{Ehrhardt2014c}. We will point out the similarities and differences in more detail in \S\ref{SEC:METHODS}. Moreover, in all \cite{Ehrhardt2014c,Ehrhardt2014,Ehrhardt2015} the parallel level sets functional has been smoothed and the problem has been solved using gradient based optimization. In contrast, here we consider the non-smooth convex formulation and propose a convex optimization algorithm for its solution.
\section{Problem Setting and Notation} 
Our derivation is carried out in a fully discrete setting where the object of interest $\MRIimage \in \RN^\MRIdimImage \subset \R^\MRIdimImage$ is sampled from a planar / volumetric MRI image. We will use this notation independently of the contrast, \ie~$\MRIimage$ might represent a $T_1$ or $T_2$ weighted image. Moreover, we follow a standard assumption for many acquisition sequences of no or negligibly small phase so that we are effectively dealing with real-valued images. An extension to complex-valued images could be done by means of a non-linear forward operator \cite{Fessler2004,Zibetti2010,Ehrhardt2015a} but this is out of scope of the present paper. Without phase, it is natural to assume that the MRI image $\MRIimage$ is non-negative which we will incorporate into the reconstruction. With the common assumption of additive Gaussian noise \cite{McVeigh1985,Macovski1996} a maximum a posteriori reconstruction with the prior proportional to $\exp(-\ParamReg\Reg)$, with functional $\FunFromTo{\Reg}{\R^\MRIdimImage}{\R}$ to be defined later, is equivalent to the minimization problem
\begin{align}
  \Argmin{\MRIimage \in \RN^\MRIdimImage}{\Bigl\{\frac12 \abs{\MRIop\MRIimage - \MRIdata}^2 + \ParamReg\Reg(\MRIimage)\Bigr\}}, \label{EQU:MRI:RECON}
\end{align}
where $\FunFromTo{\MRIop}{\R^\MRIdimImage}{\C^\MRIdimData}$ is the MRI forward operator and $\MRIdata$ the acquired data. Throughout the paper we use $\abs{\vecA}^2 \De \Adj{\vecA} \vecA$ to denote the standard norm for complex vectors with $\Adj{\vecA}$ being the Hermitian (complex conjugate transpose) of $\vecA$.

\subsection{Forward Operator for Magnetic Resonance Imaging}
The forward model in MRI is commonly assumed to be the Fourier transform $\FT$ \cite{Liang2000}. As we model our image to be real-valued but the Fourier transform acts on complex images, we embed the real into the complex space by means of an operator $\FunFromTo{\AdjRe}{\R^\MRIdimImage}{\C^\MRIdimImage}, \AdjRe(x) = x + 0 \imaginary$. It is not difficult to show that $\AdjRe$ is the adjoint of the real part restriction operator $\FunFromTo{\Re}{\C^\MRIdimImage}{\R^\MRIdimImage}, \Re(\vecA + \imaginary \vecB) \De \vecA$ when we equip the complex space with the inner product 
%$\angle{z_1}{z_2}_{\C^\MRIdimImage} = \Re (\Adj{z_1} z_2)$
$\angle{\vecA}{\vecB}_{\C^\MRIdimImage} = \Re (\Adj{\vecA} \vecB)$. 
Moreover, let $\MRIsamplingSequence \in \{1,\ldots,\MRIdimImage\}^\MRIdimData$ defining a sequence of sample locations which mimics an arbitrary MRI acquisition protocol. Then we can define a general sampling operator 
\begin{align}
  \FunFromTo{\Sampling}{\C^\MRIdimImage}{\C^\MRIdimData}, \quad (\Sampling x)_\InData = x_{\MRIsamplingSequence[\InData]},
\end{align}
where for better readability we denote the $\InData$th component of the vector $\MRIsamplingSequence$ by $\MRIsamplingSequence[\InData]$. Here, we focus on the case of practical interest, $\MRIdimData \ll \MRIdimImage$, where the number of measurements is much smaller than the number of unknowns. With such defined operators, the MRI forward operator for our model can be expressed as their concatenation
\begin{align}
  \MRIop \colon \R^\MRIdimImage \overset{\AdjRe}{\rightarrow} \C^\MRIdimImage \overset{\FT}{\rightarrow} \C^\MRIdimImage \overset{\Sampling}{\rightarrow} \C^\MRIdimData.
\end{align}
Due to the embedding $\AdjRe$ and the sampling $\Sampling$ this operator is in general not invertible. 

For the reconstruction method proposed in section \S\ref{SSEC:ADMM} we need the adjoint of $\MRIop$ which is given as
\begin{align}
  \Adj{\MRIop} \colon \C^\MRIdimData \overset{\Adj{\Sampling}}{\rightarrow} \C^\MRIdimImage \overset{\FT^{-1}}{\rightarrow} \C^\MRIdimImage \overset{\Re}{\rightarrow} \R^\MRIdimImage
\end{align}
with $\FT^{-1}$ denoting the inverse Fourier transform and $\Adj{\Sampling}$ the adjoint of the sampling operator. The latter is given by $\Adj{\Sampling}(\MRIdata) \De \sum_{\InData=1}^{\MRIdimData} \MRIdata_\InData \delta_{\InData,\MRIsamplingSequence[\InData]}$ with the Kronecker delta $\delta_{\InData,\InImage} = 1$ if $\InData=\InImage$ and $0$ otherwise, see \eg~\cite{Ehrhardt2015a}.
  
\subsection{Discrete Gradient}
The functional $\Reg$ in \eqref{EQU:MRI:RECON} encodes the a priori information in a way such that unlikely or undesirable solutions $\MRIimage$ result in a large value $\Reg(\MRIimage)$. For images it is common to penalize changes between neighbouring pixel values which can be expressed by the discrete gradient operator.

At every location $\InImage = 1, \ldots, \MRIdimImage$ we define a discrete gradient $\Grad\MRIimage_\InImage \in \GradSpace$. In the numerical simulations, we will use forward differences in two dimensions such that $\GradSpace = \R^2$ but other choices are possible, too. In general, the discrete gradient operator $\FunFromTo{\Grad}{\R^\MRIdimImage}{\GradSpace^\MRIdimImage}$ should be a linear mapping from space of images to the space of gradients. We make use of the discrete divergence operator $\FunFromTo{\Div}{\GradSpace^\MRIdimImage}{\R^\MRIdimImage}$ defined as the negative adjoint of the gradient, \ie~for all $\SubGradA \in \GradSpace^\MRIdimImage, \MRIimage \in \R^\MRIdimImage$ it holds $\angle{\Div\SubGradA}{\MRIimage}_{\R^\MRIdimImage} = \angle{\SubGradA}{- \Grad \MRIimage}_{\GradSpace^\MRIdimImage}$. For an approximation of the gradient with forward differences the matching approximation for the divergence corresponds to backward differences \cite{Aubert2001}. Moreover, let $\MatSpace \De \Lin(\GradSpace)$ be the space of linear mappings from $\GradSpace$ to $\GradSpace$. Then, we define the multiplication of a matrix-field $\MatAnisotropy \in \MatSpace^\MRIdimImage$ with a vector-field $\SubGradA \in \GradSpace^\MRIdimImage$ pointwise as $\MatAnisotropy \SubGradA \in \GradSpace^\MRIdimImage$ with $(\MatAnisotropy \SubGradA)_\InImage \De \MatAnisotropy_\InImage \SubGradA_\InImage$, a matrix-vector multiplication at the particular location.

\section{Modelling A Priori Information}\label{SEC:METHODS}
\subsection{Total Variation}
A popular regularization $\Reg$ in a variational formulation \eqref{EQU:MRI:RECON} is the total variation \cite{Rudin1992} which in our discrete setting reads
\begin{align}
  \FunFromTo{\SymbolTV}{\R^\MRIdimImage}{\R}, \quad \SymbolTV(\MRIimage) \De \sum_{\InImage=1}^\MRIdimImage \abs{\Grad \MRIimage_\InImage} \label{EQU:TV}
\end{align}
with the discrete gradient operator as defined in the previous section. The total variation has many desirable properties: it is convex and it leads to edge-preserved denoising. However, the standard formulation  does not allow to incorporate any extra a priory knowledge about the solution.

\subsection{Incorporating Structural Knowledge}
\subsubsection{A Priori Information on Location of Edges}
While the actual intensities of two MRI contrasts are very different, their structure in terms of edges is likely to be highly correlated. To incorporate the information about the location of edges extracted from one contrast, $\MRIsinfo$, into the reconstruction of the other we propose to introduce weights $\WTVweight_\InImage$ into the total variation functional.
\begin{definition}[Weighted Total Variation]~
Let $\WTVweight \in [0,1]^\MRIdimImage$ be a vector of weights. We define the \emph{weighted total variation} as
\begin{align}
  \FunFromTo{\SymbolWTV}{\R^\MRIdimImage}{\R}, \quad \SymbolWTV(\MRIimage) \De \sum_{\InImage=1}^\MRIdimImage \WTVweight_\InImage \abs{\Grad \MRIimage_\InImage}. \label{EQU:WTV}
\end{align}
\end{definition}
\begin{remark}
An option for the choice of such weights is $\WTVweight_\InImage = \ParamThresh / \abs{\Grad\MRIsinfo_\InImage}_\ParamThresh$, where $\abs{\vecA}_\ParamThresh^2 \De \abs{\vecA}^2 + \ParamThresh^2$ for some parameter $\ParamThresh > 0$. This choice results in $0 < \WTVweight_\InImage \leq 1$, with the upper bound attained when there is no side information, \ie~$\MRIsinfo = const$, hence $\abs{\Grad\MRIsinfo_\InImage} = 0$ and the lower bound approached asymptotically for $\abs{\Grad\MRIsinfo_\InImage} \rightarrow \infty$. The parameter $\ParamThresh$ controls what magnitude of an edge is considered to be \Quote{large} and what is considered to be \Quote{small}. While in general this parameter could be a spatial map, for simplicity here we assume that it is constant over space.
\end{remark}
\begin{remark}
Obviously, for the choice of uninformative weights, \ie~$\WTVweight_\InImage = 1$ for all $\InImage=1,\ldots,\MRIdimImage$, the weighted total variation functional $\SymbolWTV$ reduces to the standard total variation \eqref{EQU:TV}. Furthermore, $0 \leq \WTVweight_\InImage \leq 1$ implies that for all $\MRIimage \in \R^\MRIdimImage$ it holds $0 \leq \SymbolWTV(\MRIimage) \leq \SymbolTV(\MRIimage)$.
\end{remark}

\subsubsection{A Priori Information on Direction of Edges}
In the weighted total variation functional \eqref{EQU:WTV} we made use of the location of the edges by means of weights depending on the modulus of the gradient of the side information. However, it is reasonable to assume that these images do not only share the location but also the direction of edges modulo their sign. The latter is necessary as the actual intensity values are independent of one another such that in one image their might be a jump \Quote{up} while in the other one there is a jump \Quote{down}. 
\begin{definition}[Directional Total Variation]~Let $\AQPLdirection \in \GradSpace^\MRIdimImage$ with $0 \leq \abs{\AQPLdirection_\InImage} \leq 1$ be a vector-field and $\ProjA_{\AQPLdirection_\InImage} \De \Id - \AQPLdirection_\InImage \Adj{\AQPLdirection_\InImage}$, \ie~$\ProjA_{\AQPLdirection_\InImage} \vecA = \vecA - \angle{\AQPLdirection}{x} \AQPLdirection$. We call
\begin{align}
  \FunFromTo{\SymbolAQPL}{\R^\MRIdimImage}{\R}, \quad \SymbolAQPL(\MRIimage) \De \sum_{\InImage=1}^\MRIdimImage \abs{\ProjA_{\AQPLdirection_\InImage}\Grad \MRIimage_\InImage} \label{EQU:DTV}
\end{align}
the \emph{directional total variation}.
\end{definition}
\begin{remark}
In this paper we choose $\AQPLdirection \in \GradSpace^\MRIdimImage, \AQPLdirection_\InImage \De \Grad \MRIsinfo_\InImage / \abs{\Grad\MRIsinfo_\InImage}_\ParamThresh$ which captures the \Quote{structure} of $\MRIsinfo$ with more degrees of freedom than in the case of weighted total variation, \cf~figure \ref{FIG:GROUNDTRUTH}. As in the previous case, we will make use of an edge parameter $\ParamThresh$ that is related to the size of an edge. Similar to \eqref{EQU:WTV}, we have $0 \leq \abs{\AQPLdirection_\InImage} < 1$ with the lower bound being attained for $\abs{\Grad\MRIsinfo_\InImage} = 0$ and the upper bound approached as $\abs{\Grad\MRIsinfo_\InImage} \rightarrow \infty$. In the limit $\abs{\AQPLdirection_\InImage} \rightarrow 1$, $\ProjA_{\AQPLdirection_\InImage}$ becomes the orthogonal projection onto the orthogonal complement of $\AQPLdirection_\InImage$. Thus, in contrast to isotropic weighting of $\abs{\Grad \MRIimage_\InImage}$ in \eqref{EQU:WTV}, in the limit \eqref{EQU:DTV} penalizes only the component of $\Grad \MRIimage_\InImage$ that is orthogonal to $\AQPLdirection_\InImage$ resulting in an anisotropic weighting.
\end{remark}
\begin{remark}
The directional total variation \eqref{EQU:DTV} for $\tilde \AQPLdirection_\InImage \propto \Grad \MRIsinfo_\InImage$ is related to the parallel level sets approach \cite{Ehrhardt2014c,Ehrhardt2015,Ehrhardt2015a}. To be more precise, it was proven in \cite{Ehrhardt2015a} that
\begin{align}
  \abs{\ProjA_{\AQPLdirection_\InImage} \Grad \MRIimage_\InImage}
  = \bigl(\abs{\Grad \MRIimage_\InImage}^2 - \angle{\Grad \MRIimage_\InImage}{\tilde \AQPLdirection_\InImage}^2 \bigr)^{1/2}, \label{EQU:PL}
\end{align}
which shows that directional total variation is a special case of asymmetric parallel level sets with a different normalization of the side information $\tilde \AQPLdirection_\InImage \De (2 - \abs{\AQPLdirection_\InImage}^2)^{1/2} \AQPLdirection_\InImage$. From \eqref{EQU:PL} it can be seen that directional total variation favours parallel level sets. Indeed, on the one hand, \eqref{EQU:PL} is minimal if and only if $\Grad \MRIimage_\InImage$ is parallel to (in the span of) $\tilde \AQPLdirection_\InImage$ and hence parallel to $\Grad \MRIsinfo_\InImage$. On the other hand, as gradients are orthogonal to level sets, parallel gradients imply parallel level sets.
\end{remark}

\subsubsection{General Framework}
Both functionals \eqref{EQU:WTV} and \eqref{EQU:DTV} can be uniformly written as
\begin{align}
  \Reg(\MRIimage) = \sum_{\InImage=1}^\MRIdimImage \abs{\MatAnisotropy_\InImage \Grad \MRIimage_\InImage} ,\label{EQU:TV:STRUCTURE}
\end{align}
where the matrix-field $\MatAnisotropy \in \MatSpace^\MRIdimImage$ depends on the structural knowledge derived from the image $\MRIsinfo$. In the case of weighted total variation, 
\begin{align}
  \MatAnisotropy_\InImage = \WTVweight_\InImage \Id,
\end{align}
the matrix-field is \emph{isotropic}, \ie~it is directionally independent. On the other hand, for directional total variation, 
\begin{align}
  \MatAnisotropy_\InImage = \Id - \AQPLdirection_\InImage \Adj{\AQPLdirection_\InImage},
\end{align}
the matrix field is \emph{anisotropic} as it has principle directions along and orthogonal to $\AQPLdirection_\InImage$. As $\AQPLdirection_\InImage$ was defined to be the normalized gradient field of $\MRIsinfo$ these directions correspond to the normal and tangential direction to the level sets of $\MRIsinfo$.
\section{Algorithmic Approach}
In order to numerically solve problem \eqref{EQU:MRI:RECON} we will reformulate the problem such that it can be efficiently solved with the alternating direction method of multipliers (ADMM), see \cite{Afonso2010,Boyd2010} and references therein. As we model MRI images to be real-valued, it is efficient to perform two splits. Similar to total variation regularization, no closed-form proximal operator for priors of the form \eqref{EQU:TV:STRUCTURE} exists, thus we revert to a variant of fast gradient projection algorithm \cite{Beck2009}.

\subsection{Proximal Operator with Fast Gradient Projection}
\begin{algorithm}[t]
% \FigureSettings
\caption{Fast Gradient Projection Method for Structure-Guided Total Variation}\label{ALG:FGP}
\begin{algorithmic}[1]
\Require ~\linebreak 
	\begin{tabular}{p{1.5cm}l}
	$ \ParamReg \geq 0$ & regularization parameter\\
	$ \FGPproxpoint \in \R^\MRIdimImage$ & proximal point \\
	$ \IterA \in \N$ & number of iterations \\ 
	$ \MatAnisotropy$ & anisotropy (default = $\id$) \\
	$ \ProjA _\ContraintSetImage$ & projection onto the set $\ContraintSetImage$ (default = $\id$)\\
	$ \SubGradA ^0$ & initial dual variable  (default = 0)
	\end{tabular}
\Ensure ~\linebreak 
	\begin{tabular}{p{1.5cm}l}
	$ \MRIimage ^\IterA$ & approximation of minimizer (primal variable)\\
	$ \SubGradA ^\IterA$ & dual variable 
	\end{tabular}
\Function{FGP\_J}{$ \ParamReg, \FGPproxpoint, \IterA, \MatAnisotropy, \ProjA_C, \SubGradA ^0$}
	\State $t^0 \gets 1, \SubGradB ^0 \gets \SubGradA ^0$ \Comment{initialization}
	\For{$\iterA = 1 : \IterA $}
	\State $g ^{\iterA} \gets \ParamReg \MatAnisotropy\Grad \ProjA_\ContraintSetImage(\FGPproxpoint + \ParamReg \Div \Adj{\MatAnisotropy} \SubGradB ^{\iterA-1})$ \Comment{compute gradient \eqref{EQU:GRADDUAL}}
	\State $\SubGradA ^\iterA \gets \ProjA_\UnitBall (\SubGradB ^{\iterA-1} + \FGPstepsize g ^{\iterA})$ \Comment{update dual variable}
	\State $t ^\iterA \gets \tfrac12 \Bigl(1 + \sqrt{1 + 4 (t^{\iterA-1})^2} \Bigr)$ %$t ^\iterA \gets \tfrac12\{1 + [1 + 4 (t^{\iterA-1})^2]^{1/2}\}$ 
	\Comment{update step size}	
	\State $\SubGradB ^\iterA \gets \SubGradA ^\iterA + \frac{t ^{\iterA-1} - 1}{t ^\iterA} (\SubGradA ^\iterA - \SubGradA ^{\iterA-1})$ \Comment{Nesterov two step update}
	\EndFor
\State $\MRIimage ^\iterA \gets \ProjA_\ContraintSetImage(\FGPproxpoint + \ParamReg \Div \Adj{\MatAnisotropy} \SubGradA ^\iterA)$ \Comment{calculate final primal variable}
\State \textbf{return} $(\MRIimage ^\IterA,\SubGradA ^\IterA)$
\EndFunction
\end{algorithmic}
\end{algorithm}
Evaluation of proximal operator for structural total variation \eqref{EQU:TV:STRUCTURE}, entails solution of the following convex minimization problem
\begin{align}
	\Prox{\ParamReg\Reg + \chi_\ContraintSetImage}{\FGPproxpoint} \De \Argmin{\MRIimage\in\ContraintSetImage}{\Bigl\{\frac12 \abs{\MRIimage - \FGPproxpoint}^2 + \ParamReg\Reg(\MRIimage)\Bigr\}} \label{EQU:STV:PROX}
\end{align}
with the non-empty, closed and convex constraint set $\ContraintSetImage \subset \R^\MRIdimImage$. Although we are primarily interested in non-negativity constraints, \ie~$\ContraintSetImage = \RN^\MRIdimImage$, there is no need to be too specific at this point. Analogously to the case for usual total variation \cite{Beck2009}, structural total variation can be dualized as
\begin{align}
	\Reg(\MRIimage) = \sum_{\InImage=1}^\MRIdimImage \abs{\MatAnisotropy_\InImage \Grad \MRIimage_\InImage} = \sup_{\SubGradA\in\UnitBall} \angle{- \Div \Adj{\MatAnisotropy} \SubGradA}{\MRIimage}, \label{EQU:STV:DUAL}
\end{align}
where the supremum is taken over the unit ball in the gradient space $\UnitBall \De \{\vecA \in \GradSpace^\MRIdimImage : \abs{\vecA_\InImage} \leq 1 \}$. Substituting \eqref{EQU:STV:DUAL} into \eqref{EQU:STV:PROX} and exchanging the order of the minimum and supremum as the function is convex in $\MRIimage$ and concave in $\SubGradA$ (see \eg~\cite{Rockafellar1970}, Corollary 37.3.2) we obtain
\begin{align}
	\underset{\MRIimage \in \ContraintSetImage}{\operatorname{\phantom{g}\hspace*{-1.5mm}min}} \Bigl\{\frac12 \abs{\MRIimage - \FGPproxpoint}^2 + \ParamReg\Reg(\MRIimage)\Bigr\}
	&= \sup_{\SubGradA \in \UnitBall} \underset{\MRIimage \in \ContraintSetImage}{\operatorname{\phantom{g}\hspace*{-1.5mm}min}} \Bigl\{ \frac12 \abs{\MRIimage - \FGPproxpoint}^2 + \ParamReg \angle{-\Div\Adj{\MatAnisotropy} \SubGradA}{\MRIimage}\Bigr\} \label{EQU:SUPMIN}\\
	&= \sup_{\SubGradA \in \UnitBall} \Bigl\{ \frac12 \abs{\Minimizer{\MRIimage}(\SubGradA) - \FGPproxpoint}^2 + \ParamReg \angle{-\Div\Adj{\MatAnisotropy} \SubGradA}{\Minimizer{\MRIimage}(\SubGradA)}\Bigr\}, \label{EQU:DUALPROB}
\end{align}
where the inner minimization in \eqref{EQU:SUPMIN} has the solution $\Minimizer{\MRIimage}(\SubGradA) = \ProjA_\ContraintSetImage(h_p), h_p \De \FGPproxpoint + \ParamReg \Div \Adj{\MatAnisotropy} \SubGradA$. 
Following \cite{Beck2009}, the function under supremization in \eqref{EQU:DUALPROB} can be equivalently rewritten as
\begin{align}
	\frac12 \abs{\Minimizer{\MRIimage}(\SubGradA) - \FGPproxpoint}^2 + \ParamReg \angle{-\Div\Adj{\MatAnisotropy} \SubGradA}{\Minimizer{\MRIimage}(\SubGradA)}
	= \frac12 \abs{ h_p - \ProjA_\ContraintSetImage(h_p) }^2 - \frac12 \abs{h_p}^2 + \frac12 \abs{\FGPproxpoint}^2
\end{align}
and its gradient with respect to $\SubGradA$ is given by
\begin{align}
&- \ParamReg \MatAnisotropy \Grad [h_p - \ProjA_\ContraintSetImage(h_p)] + \ParamReg \MatAnisotropy \Grad h_p 
= \ParamReg \MatAnisotropy \Grad \ProjA_\ContraintSetImage (h_p)
= \ParamReg \MatAnisotropy \Grad \ProjA_\ContraintSetImage(\FGPproxpoint + \ParamReg \Div \Adj{\MatAnisotropy} \SubGradA). \label{EQU:GRADDUAL}
\end{align}

A variant of the fast projected gradient algorithm (with Nesterov acceleration) for solution of \eqref{EQU:DUALPROB} and hence \eqref{EQU:STV:PROX} is outlined in algorithm \ref{ALG:FGP}, where the orthogonal projection onto $\UnitBall$ is given by
\begin{align}
	\ProjA_\UnitBall(\SubGradA) = \SubGradA/\max(1,\abs{\SubGradA}).
\end{align}
\begin{remark}
As an instance of fast iterative soft thresholding algorithm (FISTA), algorithm \ref{ALG:FGP} with a step size $\FGPstepsize \leq (\ParamReg^2 \norm{\MatAnisotropy}^2 \norm{\Grad}^2)^{-1}$ converges in objective function values with rate $\mathcal O(1/k^2)$ \cite{Beck2009a} as $\ParamReg^2 \norm{\MatAnisotropy}^2 \norm{\Grad}^2$ is an upper bound on the Lipschitz constant of the gradient of the dual problem \eqref{EQU:GRADDUAL}, see \cite{Beck2009} for details. For both regularizers in this paper it holds $\norm{\MatAnisotropy} \leq 1$. Moreover, we approximate the gradient with forward and the divergence with backward differences for which we have $\norm{\Grad} \leq \sqrt 8$ in 2D space \cite{Beck2009} such that in both cases of interest $\FGPstepsize =(8 \ParamReg^2)^{-1}$ is sufficient to guarantee convergence in function values.
\end{remark}

\subsection{Double-Split Alternating Direction Method of Multipliers}\label{SSEC:ADMM}
Recall that we want to solve \eqref{EQU:MRI:RECON} 
\begin{align}
\nonumber \Argmin{\MRIimage \in \RN^\MRIdimImage}{\frac12 \abs{\MRIop \MRIimage - \MRIdata}^2 + \ParamReg \Reg(\MRIimage)}
\end{align}
with $\Reg$ as in \eqref{EQU:TV:STRUCTURE}. To fully exploit the structure of our forward operator $\MRIop = \tilde\FT \circ \Sampling, \tilde\FT \De \AdjRe \circ \FT$, we recast the problem as a constraint optimization problem 
\begin{align}
	\Minimizer{\MRIimage} \in \Argmin{\MRIimage \in \RN^\MRIdimImage}{\frac12 \abs{\Sampling \ADMMvarX - \MRIdata}^2 
	+ \ParamReg \Reg(\MRIimage)} 
	\quad \text{\st} \quad \ADMMvarX = \tilde\FT \ADMMvarZ, \quad \MRIimage = \ADMMvarZ
\end{align}
with the associated augmented Lagrangian
\begin{align}
	\ADMMlagrange(\MRIimage,\ADMMvarX,\ADMMvarZ)
	\De \frac12\abs{\Sampling \ADMMvarX - \MRIdata}^2 + \ParamReg\Reg(\MRIimage) + \frac \ADMMrho 2 \Bigl\{\abs{\ADMMvarX - \tilde\FT \ADMMvarZ + \ADMMlagmultX}^2 + \abs{\MRIimage - \ADMMvarZ + \ADMMlagmultU}^2 - \abs{\ADMMlagmultX}^2 - \abs{\ADMMlagmultU}^2 \Bigr\}
\end{align}
and $\ADMMlagmultX \in \C^\MRIdimImage,\ADMMlagmultU \in \R^\MRIdimImage$ are the scaled Lagrange multipliers. In order to make the algorithm as efficient as possible, $\MRIimage$ and $\ADMMvarX$ are associated with the first and $\ADMMvarZ$ with the second block of ADMM \cite{Afonso2010,Boyd2010}. Thus in every iteration we need to solve
\begin{align}
	\Argmin{(\MRIimage,\ADMMvarX) \in \RN^\MRIdimImage \times \C^\MRIdimImage}{\ADMMlagrange(\MRIimage,\ADMMvarX,\ADMMvarZ)} 
	\quad \text{ and } \quad
	\Argmin{\ADMMvarZ \in \R^\MRIdimImage}{\ADMMlagrange(\MRIimage,\ADMMvarX,\ADMMvarZ)}.
\end{align}
\begin{algorithm}[t]
% \FigureSettings
\newcommand{\Aindent}{\hspace*{\algorithmicindent}}
\caption{ADMM for MRI reconstruction}\label{ALG:ADMM}
\begin{algorithmic}[1]
\Require ~\linebreak
	\begin{tabular}{p{1.5cm}l}
	$ \MRIdata \in \C^\MRIdimData$ & MRI data\\
	$ \ParamReg \geq 0$ & regularization parameter\\
	$\Sampling$ & sampling \\
	$ \IterA \in \N$ & number of iterations
	\end{tabular}
\Ensure ~\linebreak
	\begin{tabular}{p{1.5cm}l}
	$\MRIimage^\IterA$ & approximate minimizer
	\end{tabular}	 
\Function{ADMM\_MRI}{$\MRIdata, \ParamReg, \Sampling, \IterA$}
	\State $ \ADMMrho \gets 1, \ADMMvarZ^0, \ADMMlagmultX^0, \ADMMlagmultU^0 \gets 0$ \Comment{initialize variables}
	\For{$ \iterA = 1 : \IterA$}\hfill ~\linebreak
	  \Aindent\Aindent update first block
		\State \Aindent $\MRIimage^\iterA \gets \Prox{\ParamReg/\ADMMrho \Reg + \chi_{\RN^\MRIdimImage}}{\ADMMvarZ^{\iterA-1} 
		  - \ADMMlagmultU^{\iterA-1}}$ \Comment{apply algorithm \ref{ALG:FGP}}
		\State \Aindent $\ADMMvarX^\iterA \gets (\Adj{\Sampling}\Sampling + \ADMMrho \Id)^{-1} 
		  [\Adj{\Sampling} \MRIdata + \ADMMrho (\FT \ADMMvarZ^{\iterA-1} - \ADMMlagmultX^{\iterA-1})]$ \Comment{component-wise scaling}\linebreak
		\Aindent\Aindent update second block
		\State \Aindent $\ADMMvarZ^\iterA \gets \frac12 [\Re\FT^{-1} (\ADMMvarX^\iterA + \ADMMlagmultX^{\iterA-1}) 
		  + \MRIimage^\iterA + \ADMMlagmultU^{\iterA-1}]$ \Comment{averaging step} \linebreak
		\Aindent\Aindent update Lagrange multipliers
		\State \Aindent $\ADMMlagmultX^\iterA \gets \ADMMlagmultX^{\iterA-1} 
		  + \ADMMvarX^\iterA - \FT \ADMMvarZ^\iterA$
		\State \Aindent $\ADMMlagmultU^\iterA \gets \ADMMlagmultU^{\iterA-1} 
		  + \MRIimage^\iterA - \ADMMvarZ^\iterA$\hfill ~\linebreak
		\Aindent\Aindent update $ \ADMMrho $ according to \cite{Boyd2010}
	\EndFor
\State \textbf{return} $\MRIimage^\IterA$
\EndFunction
\end{algorithmic}
\end{algorithm}
As the first minimization problem decouples in $\MRIimage$ and $\ADMMvarX$, we obtain three update steps for ADMM, the first two of which can be performed in parallel
\begin{align}
	&\Argmin{\MRIimage \in \RN^\MRIdimImage}{\Bigl\{\frac \ADMMrho 2 \abs{\MRIimage - \ADMMvarZ + \ADMMlagmultU}^2 + \ParamReg\Reg(\MRIimage)\Bigr\}}
	= \Prox{\alpha/\ADMMrho\Reg + \chi_{\RN^\MRIdimImage}}{\ADMMvarZ - \ADMMlagmultU}\\
	&\Argmin{\ADMMvarX \in \C^\MRIdimImage}{\Bigl\{\frac12\abs{\Sampling \ADMMvarX - \MRIdata}^2 + \frac\ADMMrho2 \abs{\ADMMvarX - \tilde\FT \ADMMvarZ + \ADMMlagmultX}^2}\Bigr\}
	= (\Adj{\Sampling}\Sampling + \ADMMrho \Id)^{-1} [\Adj{\Sampling} \MRIdata + \ADMMrho (\FT \ADMMvarZ - \ADMMlagmultX) ] \label{EQU:ADMM:UPDATE:X}\\
	&\Argmin{\ADMMvarZ \in \R^\MRIdimImage}{\Bigl\{\frac12 \abs{\ADMMvarX - \tilde\FT \ADMMvarZ + \ADMMlagmultX}^2 + \frac12\abs{\MRIimage - \ADMMvarZ + \ADMMlagmultU}^2}\Bigr\}
	= \frac12 [\Adj{\tilde\FT} (\ADMMvarX + \ADMMlagmultX) + \MRIimage + \ADMMlagmultU].
\end{align}
In \eqref{EQU:ADMM:UPDATE:X} we used that $\tilde\FT \ADMMvarZ = \FT \ADMMvarZ$ for real $\ADMMvarZ$. It should be noted that both $\Adj{\Sampling}\Sampling$ and $\ADMMrho \Id$ are diagonal matrices, so that the matrix inversion in \eqref{EQU:ADMM:UPDATE:X}, $(\Adj{\Sampling}\Sampling +\ADMMrho \Id)^{-1}$, reduces to a component-wise scaling and is therefore computationally efficient. The final ADMM algorithm can be found in \ref{ALG:ADMM}. In each iteration of the algorithm, we apply once the discrete Fourier transform and its inverse as well as the proximal operator via algorithm \ref{ALG:FGP}. After each iteration, if the primal and dual residual are too far apart, we adjust the parameter $\ADMMrho$ according to the guidelines in \cite{Boyd2010}.
\begin{remark}
  In vector notation, the double split can be written as $(\MRIimage; \ADMMvarX) = \mathcal G\ADMMvarZ$ where $\mathcal G \De (\Id;\tilde\FT)$ has full column rank. If the we compute the proximal operator with sufficient accuracy, \ie~the errors are absolutely summable, and $\ADMMrho$ is constant, then algorithm \ref{ALG:ADMM} converges to a solution of \eqref{EQU:MRI:RECON} \cite{Eckstein1992}. Numerically, we observe convergence for both $\MRIimage$ and $\ADMMrho$.
\end{remark}
\section{Numerical Experiments}
\subsection{Technical Details}
\subsubsection{Data and Algorithms}
We numerically test the two extensions for total variation to incorporate structural information with six datasets that are either based on the Shepp-Logan phantom, realistically simulated MRI from BrainWeb \cite{Cocosco1997} and clinical MRI images from a patient, \cf~figure \ref{FIG:GROUNDTRUTH}. We simulate the MRI data by sampling from the discrete Fourier transform in a variety of ways including Cartesian sampling (equidistantly and randomly undersampled), radial sampling (equidistantly spaced radial spokes, golden angle \cite{Winkelmann2007}) and spiral sampling (variable density and phyllotaxis \cite{Vogel1979}). In all cases we added Gaussian noise to the complex-valued MRI data with standard deviation scaled such that for fully sampled data the expected $\ell^2$-norm of the noise is 5\% of the $\ell^2$-norm of the noise-free data. 

Both algorithms \ref{ALG:FGP} and \ref{ALG:ADMM} have been implemented in MATLAB R2015a. The algorithms and the datasets used in this paper are available as supplementary material.
\begin{figure}%
\FigureSettings
\Portable{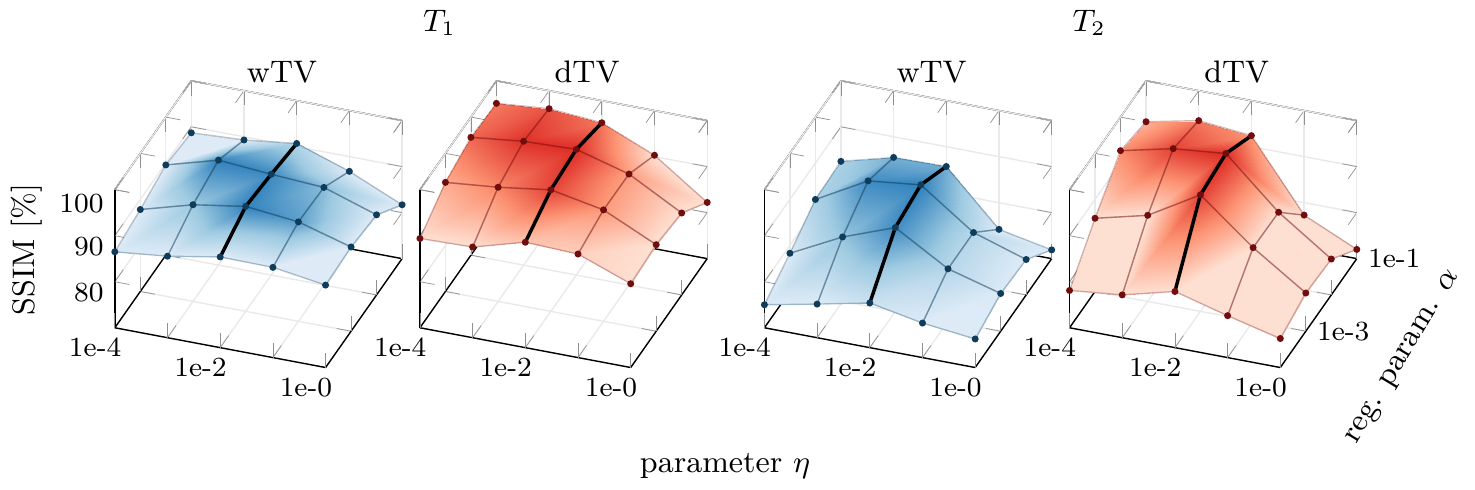}{
\newcommand{\File}[2]{pics/results/MRI2_Jun15_BrainWebNormalTransverse1_20150724_1642/2015-07-28/stats_#1_#2.txt}
\begin{tikzpicture}%
\setlength{\Width}{3.1cm}
\setlength{\Height}{3.2cm}

\pgfplotsset{
  grid=major, grid style={gray!20}, width=4.5cm, height=4.5cm, view={20}{40}, 
  xtick = {1,2,3,4,5}, xticklabels = {1e-4, , 1e-2, , 1e-0}, xticklabel style={anchor=north, xshift=-.2cm}, xlabel = {}, ytick = {1,2,3,4}, yticklabels = {}, yticklabel style={anchor=west, xshift=0cm},
  zmin = 70, zmax = 100, ztick = {100,90,80}, every tick label/.append style={font=\scriptsize}}  
\newcommand{\PlotSurf}[4]{
  \addplot3 [surf, shader=interp, colormap name={#4}] table [x = ieta, y = ialpha, z=#2] {#1};
  \addplot3 [mesh, color=#3, opacity=.3] table [x = ieta, y = ialpha, z=#2] {#1};
  \addplot3 [only marks, mark size=.8pt, color=#3] table [x = ieta, y = ialpha, z=#2] {#1}; 
  \addplot3 [color=black,line width=1pt, restrict x to domain={3:3}] table [x = ieta, y = ialpha, z=#2] {#1};}
  
  \FirstPic
  \DrawVText{\x-9mm}{\y+12mm}{SSIM [\%]} 
  \DrawText{\x+33mm}{\y+35mm}{$T_1$} 
  \DrawText{\x+17mm}{\y+30mm}{$ \SymbolWTV $} 
  \begin{axis}[at={(\x,\y)},point meta min=90,point meta max=92]
  \PlotSurf{\File{WTV}{T1}}{ssim}{\ColourWTV!50!black}{Blues-3}
  \end{axis}
  
  \NextPic
  \DrawText{\x+17mm}{\y+30mm}{$ \SymbolAQPL $} 
  \begin{axis}[at={(\x,\y)},zticklabels = {},point meta min=92,point meta max=97]
  \PlotSurf{\File{AQPL}{T1}}{ssim}{\ColourAQPL!50!black}{Reds-3}
  \end{axis}
  
  \addtolength{\x}{4mm}
  
  \NextPic
  \DrawText{\x+33mm}{\y+35mm}{$T_2$} 
  \DrawText{\x+17mm}{\y+30mm}{$ \SymbolWTV $} 
  \begin{axis}[at={(\x,\y)},zticklabels = {},point meta min=78,point meta max=89]
  \PlotSurf{\File{WTV}{T2}}{ssim}{\ColourWTV!50!black}{Blues-3}
  \end{axis}
  
  \NextPic
  \DrawText{\x+17mm}{\y+30mm}{$ \SymbolAQPL $} 
  \begin{axis}[at={(\x,\y)},y label style={yshift=.4cm,xshift=.4cm,rotate=58},yticklabels = {, 1e-3, , 1e-1},
  zticklabels = {},point meta min=89,point meta max=96]
  \PlotSurf{\File{AQPL}{T2}}{ssim}{\ColourAQPL!50!black}{Reds-3}
  \end{axis}
  
  \DrawText{\x-35mm}{\y-10mm}{parameter $ \ParamThresh $}
  \DrawText{\x+34mm}{\y+1mm}{\rotatebox{59}{reg. param. $ \ParamReg $}}
\end{tikzpicture}}
\caption{Reconstruction quality in terms of SSIM using $\SymbolWTV$ and $\SymbolAQPL$ in dependence of the parameter $\ParamThresh$ and the regularization parameter $\ParamReg$ for the data set \NameBrainWebA~and radial sampling. In all cases $\ParamThresh$ = 1e-2 yields the best results in terms of SSIM.} \label{FIG:RES:PARAMETER:ETA}
\end{figure}

\subsubsection{Quality Measures and Parameter Selection}
We evaluate the results in terms of the peak signal-to-noise (PSNR) and the structural similarity index (SSIM) \cite{Wang2004} both available in the image processing toolbox in MATLAB R2015a.

The regularization parameter $\ParamReg$ and the edge parameter $\ParamThresh$ are chosen to maximize the SSIM between the reconstructed result and the ground truth.

\subsection{Results}
\subsubsection{Parameter Estimation}
Both proposed extensions of total variation have a parameter $\ParamThresh$ that relates to the magnitude of the gradients in the side information. Figure \ref{FIG:RES:PARAMETER:ETA} shows the SSIM of the reconstructions of $T_1$ and $T_2$ weighted images from radially sampled \NameBrainWebA~data using both structure enhancing regularizers in dependence of the regularization parameter $\ParamReg$ and the edge parameter $\ParamThresh$. In all cases the best results are obtained for $\ParamThresh$ = 1e-2 which corresponds to approximately 1\% of the maximal gradient magnitude. For a large value of $\ParamThresh$---in this example approximately 1---both regularizers perform the same and both coincide with total variation (not shown). Similar plots were obtained for the other data sets and sampling patterns and hence will be omitted for brevity. In what follows the edge parameter $\ParamThresh$ is always chosen to be 1e-2.

\subsubsection{Visual Assessment}
{
\newcommand{\FolderDataLocal}{\FolderDataA}
\newcommand{\DrawImLocal}{\DrawImA}
\newcommand{\DrawZoomLocal}{\DrawZoomA}
\newcommand{\DataSetName}{\DataSetA}
\begin{figure}%
\FigureSettings
{\newcommand{\SamplingName}{cartesianX_random_0_8}
\Portable{exp_fig_results_SheppLogan_\SamplingName}{\DrawResults}}
\caption{Reconstructions and their close-ups of \NameSheppLogan~$T_1$~(top two rows) and $T_2$ (bottom two rows) for sampling pattern shown in top left corner. While both $\SymbolWTV$ and $\SymbolAQPL$ visually improve the reconstructions compared to no prior and $\SymbolTV$, the result for $\SymbolAQPL$ shows significantly fewer artefacts.}\label{FIG:RES:SHEPPLOGAN}
\vspace*{3mm}
{\newcommand{\SamplingName}{spiral2}
\Portable{exp_fig_results_SheppLogan_\SamplingName}{\DrawResultsToneSSIM}}
\caption{Reconstructions of \NameSheppLogan~$T_1$ with spiral sampling (top row) and their SSIM maps (bottom row). The SSIM index is the mean of the respective SSIM map. As in figure \ref{FIG:RES:SHEPPLOGAN}, both $ \SymbolWTV $ and $ \SymbolAQPL$ improve on no prior and $ \SymbolTV $.}\label{FIG:RES:SHEPPLOGAN:SSIM}
\vspace*{3mm}
{\newcommand{\SamplingName}{spiralPhyll2}
\Portable{exp_fig_results_SheppLogan_\SamplingName}{\DrawResultsToneError}}
\caption{Reconstructions of \NameSheppLogan~$T_1$ (top row) for sampling along the spiral phyllotaxis. It can be seen from the difference images (bottom row) that incorporating structure from $T_2$ greatly improves the reconstruction. The PSNR is proportional to the logarithm of the $\ell^2$-norm of the respective difference image.} 
\label{FIG:RES:SHEPPLOGAN:ERRORS}%
\end{figure}%
}
{
\newcommand{\FolderDataLocal}{\FolderDataB}
\newcommand{\DrawImLocal}{\DrawImB}
\newcommand{\DrawZoomLocal}{\DrawZoomB}
\newcommand{\DataSetName}{\DataSetB}
\begin{figure}%
\FigureSettings
{\newcommand{\SamplingName}{radial_golden16}
\Portable{exp_fig_results_BrainWebTransverse1_\SamplingName}{\DrawResults}}
\caption{Reconstructions for \NameBrainWebA~$T_1$ (top two rows) and $T_2$ (bottom two rows) for golden-angle radial sampling show that utilizing directional information from the other contrast significantly improves the reconstructions on a high resolution level.} \label{FIG:RES:BRAINWEBTRANS1}%
\vspace*{3mm}
{\newcommand{\SamplingName}{cartesianX_random_0_8}
\Portable{exp_fig_results_BrainWebTransverse1_\SamplingName}{\DrawResultsToneSSIM}}
\caption{The reconstructed images of \NameBrainWebA~$T_1$ and their corresponding SSIM maps demonstrate that by incorporating more a priori knowledge (from left to right), the artefacts from random Cartesian sampling get reduced.} \label{FIG:RES:BRAINWEBTRANS1:SSIM}%
\vspace*{3mm}
{\newcommand{\SamplingName}{spiral2}
\Portable{exp_fig_results_BrainWebTransverse1_\SamplingName}{\DrawResultsTtwoError}
\caption{Reconstructions of \NameBrainWebA~$T_2$ (top row) for spiral sampling and the difference images (bottom row).} \label{FIG:RES:BRAINWEBTRANS1:ERRORS}}%
\end{figure}%
}
\begin{figure}%
\FigureSettings
{\newcommand{\FolderDataLocal}{\FolderDataC}
\newcommand{\DrawImLocal}{\DrawImC}
\newcommand{\DrawZoomLocal}{\DrawZoomC}
\newcommand{\DataSetName}{\DataSetC}
\newcommand{\SamplingName}{cartesianX7}
\Portable{exp_fig_results_BrainWebTransverse2_\SamplingName}{\DrawResults}
\caption{Combining the contrasts leads to noticeably improved reconstructions for \NameBrainWebB~with equidistant Cartesian undersampling where only every seventh line has been sampled.} \label{FIG:RES:BRAINWEBTRANS2}}%
\vspace*{3mm}
{\newcommand{\FolderDataLocal}{\FolderDataD}
\newcommand{\DrawImLocal}{\DrawImD}
\newcommand{\DrawZoomLocal}{\DrawZoomD}
\newcommand{\DataSetName}{\DataSetD}
\newcommand{\SamplingName}{cartesianX_random_0_8}
\Portable{exp_fig_results_BrainWebCoronal_\SamplingName}{\DrawResults}}
\caption{Reconstructions for data set \NameBrainWebC~with random Cartesian sampling. While exploiting structural information from the other contrast already greatly enhances the image quality, it can be seen from the close-ups that using directional information allows to recover much greater level of detail. \label{FIG:RES:BRAINWEBCORONAL}}
\end{figure}%
{
\newcommand{\FolderDataLocal}{\FolderDataE}
\newcommand{\DrawImLocal}{\DrawImE}
\newcommand{\DataSetName}{\DataSetE}
\begin{figure}%
\FigureSettings
{\newcommand{\SamplingName}{cartesianX_random_0_32}
\newcommand{\DrawZoomLocal}{\DrawZoomEA}
\Portable{exp_fig_results_NinonTransverse_\SamplingName}{\DrawResults}}
\caption{Results for \NameNinonA~with random Cartesian sampling.} \label{FIG:RES:NINONTRANSVERSE:A}%
\vspace*{3mm}
{\newcommand{\SamplingName}{cartesianY4}
\newcommand{\DrawZoomLocal}{\DrawZoomEB}
\Portable{exp_fig_results_NinonTransverse_\SamplingName}{\DrawResultsTtwo}
\caption{Results for $T_2$-weighted \NameNinonA~where every fourth line was sampled.} \label{FIG:RES:NINONTRANSVERSE:B}}%
\vspace*{3mm}
{\newcommand{\SamplingName}{spiral8}
\newcommand{\DrawZoomLocal}{\DrawZoomEC}
\Portable{exp_fig_results_NinonTransverse_\SamplingName}{\DrawResultsToneSSIM}
\caption{Reconstructions of \NameNinonA~$T_1$ with spiral sampling (top row) and the corresponding SSIM maps (bottom row).} \label{FIG:RES:NINONTRANSVERSE:C}}%
\end{figure}%
}
\begin{figure}%
\FigureSettings
\newcommand{\FolderDataLocal}{\FolderDataF}
\newcommand{\DrawImLocal}{\DrawImF}
\newcommand{\DrawZoomLocal}{\DrawZoomF}
\newcommand{\DataSetName}{\DataSetF}
\newcommand{\SamplingName}{cartesianX_random_0_32}
\Portable{exp_fig_results_NinonSagital_\SamplingName}{\DrawResults}
\caption{Results for \NameNinonB~with random Cartesian sampling.} \label{FIG:RES:NINONSAGITAL}%
\end{figure}%
Figures \ref{FIG:RES:SHEPPLOGAN} - \ref{FIG:RES:NINONSAGITAL} show results of reconstructions of $T_1$ and $T_2$ weighted images of the six ground truth image sets depicted in figure \ref{FIG:GROUNDTRUTH} using different sampling schemes. Whenever appropriate we include close-ups or SSIM maps and difference images to aid quantitative comparison. While most of the images speak for themselves and some observations are included in the captions we would like to make some general comments. In all of the aforementioned figures, but probably most visibly in figure \ref{FIG:RES:BRAINWEBTRANS2}, incorporating structural knowledge from the other contrast does visually improve the reconstruction using either $\SymbolWTV$ or $\SymbolAQPL$. When comparing $\SymbolWTV$ and $\SymbolAQPL$, one notices that while $\SymbolWTV$ results in patchy images, $\SymbolAQPL$ is able to recover smooth structures accurately. Moreover, including the directional information yields another level of improvement of fine details.

\subsubsection{Quantitative Assessment}
\begin{figure}
\FigureSettings
{\newcommand{\best}{\bf}
\begin{center}
\captionof{table}{Quantitative analysis of the results where the statistics are taken over all test cases and the best result is printed bold. On average, $\SymbolAQPL$ outperforms $\SymbolTV$ by around  6dB in PSNR and 8\% in SSIM for both contrasts.} \label{TAB:RESULTS}
\begin{tabular}{crcccccccc}
\toprule
\input{stats_best_edited.txt}
\bottomrule
\end{tabular}\\[2mm]
\end{center}}
\vspace*{3mm}
{\Portable{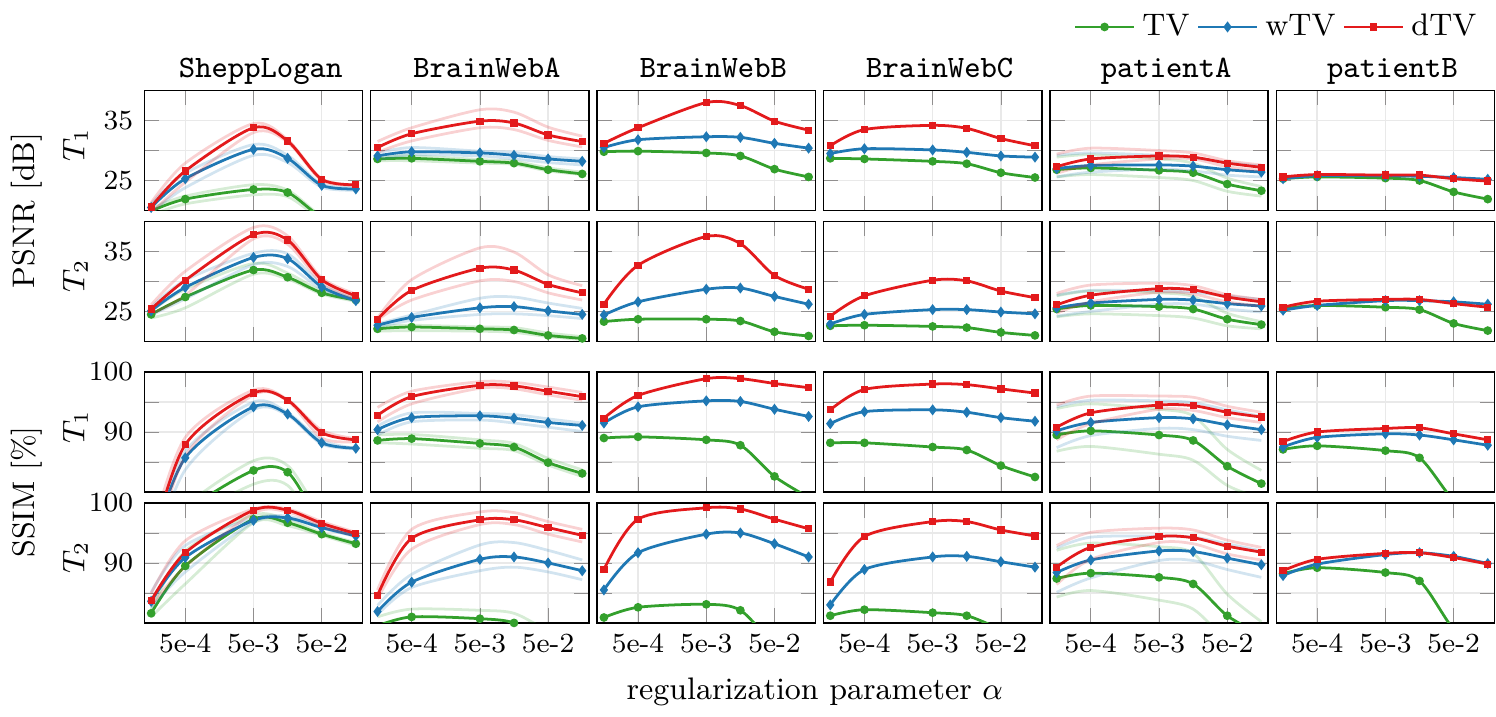}{
\newcommand{\DrawPlot}[7]{
  \addplot[color=#3, line width=.8pt, smooth, mark=#4, mark size=#6, mark options={color=#3,fill=#3,line width=0.5pt,rotate=#5},opacity=#7] table [x=No, y=#2] {pics/results/#1.txt};}%

\newcommand{\DrawPlotGroup}[6]{
  \DrawPlot{#1}{#2_MIN}{#3}{#4}{#5}{0pt}{.2}
  \DrawPlot{#1}{#2_MAX}{#3}{#4}{#5}{0pt}{.2}
  \DrawPlot{#1}{#2_MEA}{#3}{#4}{#5}{#6}{1}}

\newcommand{\DrawPlots}[1]{
  \DrawPlotGroup{stats_#1_\CONTRAST_\QUALITYMEASURE}{TV}{\ColourTV}{\MarkTV}{\MarkAngleTV}{\MarkSizeTV}
  \DrawPlotGroup{stats_#1_\CONTRAST_\QUALITYMEASURE}{WTV}{\ColourWTV}{\MarkWTV}{\MarkAngleWTV}{\MarkSizeWTV}
  \DrawPlotGroup{stats_#1_\CONTRAST_\QUALITYMEASURE}{AQPL}{\ColourAQPL}{\MarkAQPL}{\MarkAngleAQPL}{\MarkSizeAQPL}}

\newcommand{\DrawPlotsSingle}[1]{
  \DrawPlot{stats_#1_\CONTRAST_\QUALITYMEASURE}{TV_1}{\ColourTV}{\MarkTV}{\MarkAngleTV}{\MarkSizeTV}{1}
  \DrawPlot{stats_#1_\CONTRAST_\QUALITYMEASURE}{WTV_1}{\ColourWTV}{\MarkWTV}{\MarkAngleWTV}{\MarkSizeWTV}{1}
  \DrawPlot{stats_#1_\CONTRAST_\QUALITYMEASURE}{AQPL_1}{\ColourAQPL}{\MarkAQPL}{\MarkAngleAQPL}{\MarkSizeAQPL}{1}}

\pgfplotsset{every axis legend/.append style={cells={anchor=west}, at={(1,1.3)}, anchor=south east, draw=none},grid=major, grid style={gray!20}, width=3.8cm, height=2.8cm, legend columns=3,/pgfplots/xtick={2,4,6}, xmin=.8, xmax=7.2, every tick label/.append style={font=\scriptsize}}%
      
\setlength{\Width}{2.3cm}
\setlength{\Height}{1.33cm}

\newcommand{\DrawLabelsTop}[1]{
  \DrawText{\x+.5*\Width}{\y+\Height+.1cm}{\phantom{g}#1\phantom{f}}}  
\newcommand{\DrawLabelsLeft}[1]{
  \DrawVText{\x-1.2cm}{\y-0cm}{#1}}  
\newcommand{\DrawLabelsContrast}[1]{
  \DrawVText{\x-.7cm}{\y+.5*\Height}{#1}}
  
\begin{tikzpicture}%
\newcommand{\QUALITYMEASURE}{psnr} %ssim
\newcommand{\CONTRAST}{T1} %T2
\pgfplotsset{xticklabels={}}
\FirstPic\DrawLabelsLeft{PSNR [dB]}\DrawLabelsTop{\NameSheppLogan}\DrawLabelsContrast{$T_1$}
\pgfplotsset{ymin=20,ymax=40,ytick={25,30,35},yticklabels={25,,35}}
  \begin{axis}[at={(\x,\y)}]
    \DrawPlots{phantom}
  \end{axis}
\NextPic\DrawLabelsTop{\NameBrainWebA}
  \begin{axis}[at={(\x,\y)}, yticklabels={}]
    \DrawPlots{BrainWebNormalTransverse1}
  \end{axis}
\NextPic\DrawLabelsTop{\NameBrainWebB}
  \begin{axis}[at={(\x,\y)}, yticklabels={}]
    \DrawPlotsSingle{BrainWebNormalTransverse2}
  \end{axis}
\NextPic\DrawLabelsTop{\NameBrainWebC}
  \begin{axis}[at={(\x,\y)}, yticklabels={}]
    \DrawPlotsSingle{BrainWebNormalCoronal1}
  \end{axis}
\NextPic\DrawLabelsTop{\NameNinonA}
  \begin{axis}[at={(\x,\y)}, yticklabels={}]
    \DrawPlots{NinonTransverse}
  \end{axis}
\NextPic\DrawLabelsTop{\NameNinonB}
  \begin{axis}[at={(\x,\y)}, yticklabels={}]
    \DrawPlotsSingle{NinonSagital}
    \legend{\SymbolTV,\SymbolWTV,\SymbolAQPL}
  \end{axis}
\renewcommand{\CONTRAST}{T2} %T2
\NextRow\DrawLabelsContrast{$T_2$}
  \begin{axis}[at={(\x,\y)}]
    \DrawPlots{phantom}
  \end{axis}
\NextPic
  \begin{axis}[at={(\x,\y)}, yticklabels={}]
    \DrawPlots{BrainWebNormalTransverse1}
  \end{axis}
\NextPic
  \begin{axis}[at={(\x,\y)}, yticklabels={}]
    \DrawPlotsSingle{BrainWebNormalTransverse2}
  \end{axis}
\NextPic
  \begin{axis}[at={(\x,\y)}, yticklabels={}]
    \DrawPlotsSingle{BrainWebNormalCoronal1}
  \end{axis}
\NextPic
  \begin{axis}[at={(\x,\y)}, yticklabels={}]
    \DrawPlots{NinonTransverse}
  \end{axis}
\NextPic
  \begin{axis}[at={(\x,\y)}, yticklabels={}]
    \DrawPlotsSingle{NinonSagital}
  \end{axis}
\addtolength{\y}{-2mm}
\renewcommand{\QUALITYMEASURE}{ssim} %ssim
\renewcommand{\CONTRAST}{T1}
\pgfplotsset{xticklabels={},ymin=.8,ymax=1,ytick={.85,.9,.95,1},
  yticklabels={,90,,100}}
\NextRow\DrawLabelsLeft{SSIM [\%]}\DrawLabelsContrast{$T_1$}
  \begin{axis}[at={(\x,\y)}]
    \DrawPlots{phantom}
  \end{axis}
\NextPic
  \begin{axis}[at={(\x,\y)}, yticklabels={}]
    \DrawPlots{BrainWebNormalTransverse1}
  \end{axis}
\NextPic
  \begin{axis}[at={(\x,\y)}, yticklabels={}]
    \DrawPlotsSingle{BrainWebNormalTransverse2}
  \end{axis}
\NextPic
  \begin{axis}[at={(\x,\y)}, yticklabels={}]
    \DrawPlotsSingle{BrainWebNormalCoronal1}
  \end{axis}
\NextPic
  \begin{axis}[at={(\x,\y)}, yticklabels={}]
    \DrawPlots{NinonTransverse}
  \end{axis}
\NextPic
  \begin{axis}[at={(\x,\y)}, yticklabels={}]
    \DrawPlotsSingle{NinonSagital}
  \end{axis}
\renewcommand{\CONTRAST}{T2} %T2
\pgfplotsset{xticklabels={5e-4,5e-3,5e-2}}
\NextRow\DrawLabelsContrast{$T_2$}
  \begin{axis}[at={(\x,\y)}]
    \DrawPlots{phantom}
  \end{axis}
\NextPic
  \begin{axis}[at={(\x,\y)}, yticklabels={}]
    \DrawPlots{BrainWebNormalTransverse1}
  \end{axis}
\NextPic
  \begin{axis}[at={(\x,\y)}, yticklabels={}]
    \DrawPlotsSingle{BrainWebNormalTransverse2}
  \end{axis}
\NextPic
  \begin{axis}[at={(\x,\y)}, yticklabels={}]
    \DrawPlotsSingle{BrainWebNormalCoronal1}
  \end{axis}
\NextPic
  \begin{axis}[at={(\x,\y)}, yticklabels={}]
    \DrawPlots{NinonTransverse}
  \end{axis}
\NextPic
  \begin{axis}[at={(\x,\y)}, yticklabels={}]
    \DrawPlotsSingle{NinonSagital}
  \end{axis}
\DrawText{.5\linewidth-1cm}{\y - .7cm}{regularization parameter $ \ParamReg $}
\end{tikzpicture}}
\caption{Quantitative analysis of the results for the quality measures PSNR (top) and SSIM (bottom) with respect to the regularization parameter $\ParamReg$. Including structural knowledge into the reconstruction does not only improve the reconstruction for the \textit{optimal choice} of regularization parameter but it also makes it \textit{more robust} as the results are consistently better for all shown choices of regularization parameter.} \label{FIG:RES:QUANT}}
\vspace*{3mm}
{\Portable{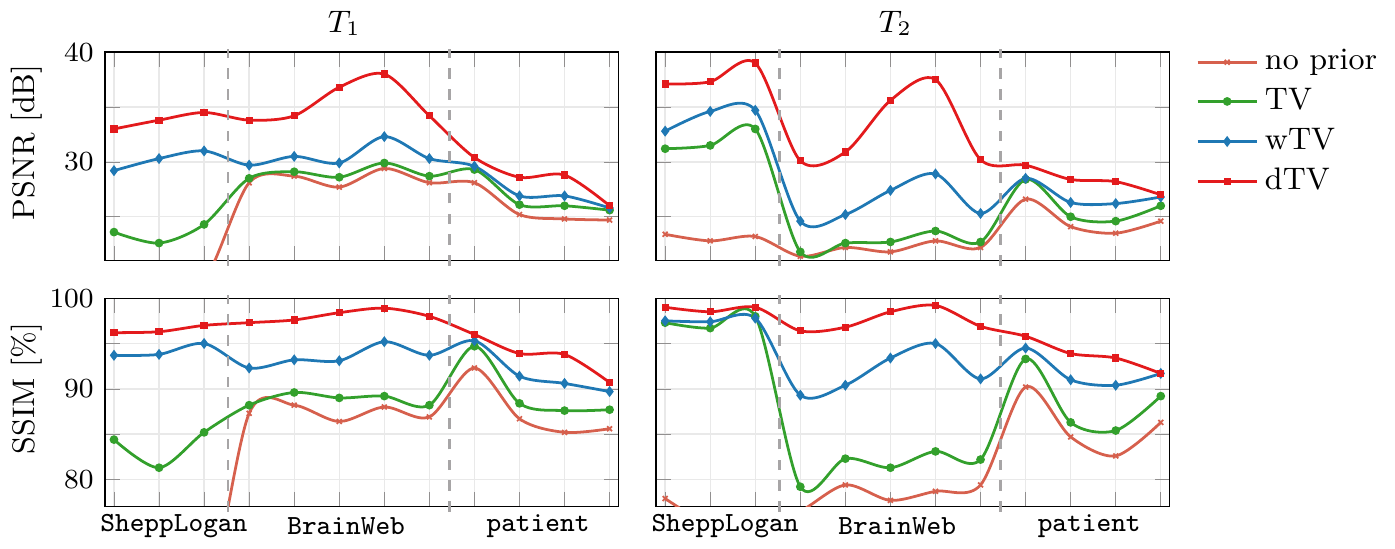}{
\newcommand{\DrawPlot}[6]{
  \addplot[color=#2, line width=.8pt, smooth, mark=#3, mark size=#5, mark options={color=#2,fill=#2,line width=0.5pt,rotate=#4},opacity=#6] table [x=No, y=#1] {pics/results/stats_summary.txt};}%

\newcommand{\DrawPlots}[2]{
  \DrawPlot{ES#1#2}{\ColourES}{\MarkES}{\MarkAngleES}{\MarkSizeES}{1}
  \DrawPlot{TV#1#2}{\ColourTV}{\MarkTV}{\MarkAngleTV}{\MarkSizeTV}{1}
  \DrawPlot{WTV#1#2}{\ColourWTV}{\MarkWTV}{\MarkAngleWTV}{\MarkSizeWTV}{1}
  \DrawPlot{AQPL#1#2}{\ColourAQPL}{\MarkAQPL}{\MarkAngleAQPL}{\MarkSizeAQPL}{1}}

\pgfplotsset{every axis legend/.append style={cells={anchor=west}, at={(1.03,1.08)}, anchor=north west, draw=none}, grid=major, grid style={gray!20}, width=6.8cm, height=3.7cm, legend columns=1,xtick={1,2,...,12},xmin=.8,xmax=12.2, every tick label/.append style={font=\scriptsize}}%

\newcommand{\DrawLabels}[2]{
  \DrawVText{\x-.8cm}{\y+1.2cm}{#1}
  \DrawText{\x+2.4cm}{\y+2.4cm}{\phantom{g}#2\phantom{f}}}
      
\newcommand{\DrawDataLabels}{
    \draw node at (\x+.7cm,\y-.2cm) {\scriptsize \NameSheppLogan};
    \draw node at (\x+2.45cm,\y-.2cm) {\scriptsize \NameBrainWeb};
    \draw node at (\x+4.4cm,\y-.2cm) {\scriptsize \NameNinon};}

\newcommand{\DrawDataLines}{
    \draw[line width=.8pt,color=gray!80,dashed] (\x+1.25cm,\y-.05cm) -- (\x+1.25cm,\y+2.15cm);
    \draw[line width=.8pt,color=gray!80,dashed] (\x+3.50cm,\y-.05cm) -- (\x+3.50cm,\y+2.15cm);}

\setlength{\Width}{5.6cm}
\setlength{\Height}{2.5cm}
\begin{tikzpicture}%
\pgfplotsset{xticklabels={},ymin=21,ymax=40,ytick={25,30,35,40},yticklabels={,30,,40}}
\FirstPic\DrawLabels{PSNR [dB]}{$T_1$}
  \begin{axis}[at={(\x,\y)}]
    \DrawPlots{T1}{psnr}
  \end{axis}
  \DrawDataLines
  
\NextPic\DrawLabels{}{$T_2$}
  \begin{axis}[at={(\x,\y)},yticklabels={}]
    \DrawPlots{T2}{psnr}
    \legend{no prior,\SymbolTV,\SymbolWTV,\SymbolAQPL}
  \end{axis}
  \DrawDataLines  
  
\pgfplotsset{xticklabels={},ymin=77,ymax=100,ytick={80,85,90,95,100},yticklabels={80,,90,,100}}
\NextRow\DrawLabels{SSIM [\%]}{}
  \begin{axis}[at={(\x,\y)}]
    \DrawPlots{T1}{ssim}
  \end{axis}
  \DrawDataLines
  \DrawDataLabels
  
\NextPic
  \begin{axis}[at={(\x,\y)},yticklabels={}]
    \DrawPlots{T2}{ssim}
  \end{axis}
  \DrawDataLines
  \DrawDataLabels
\end{tikzpicture}}
\caption{Quantitative summary of all results: While the relative performance of the methods depends strongly on the data set and the sampling scheme, incorporating structural knowledge in most cases significantly improves the results. Moreover, the best results have consistently been obtained making use of directional edge information. \label{FIG:RES:SUMMARY}}}
\end{figure}

Quantitative analysis of the results is summarized in figures \ref{FIG:RES:QUANT} and \ref{FIG:RES:SUMMARY}, and table \ref{TAB:RESULTS}. Figure \ref{FIG:RES:QUANT} shows the reconstruction quality for all six test cases in dependence of the regularization parameter. Whenever more than one sampling scheme was used, the solid line corresponds to the mean performance with the worst and best performance indicated by shaded lines. For all test cases, but especially for $T_1$-weighted \NameSheppLogan~and  \NameBrainWeb,~$\SymbolWTV$ and $\SymbolAQPL$ strongly outperforms the standard total variation. Moreover, the curves are layered which means that the results are not only better for one choice of regularization parameter but for all choices shown. 

Figure \ref{FIG:RES:SUMMARY} shows the performance for the optimal value of the regularization parameter for all test cases (data sets and sampling schemes). Also here, the curves are layered, meaning that in every test case $\SymbolAQPL$ outperforms all the other methods. The average performance can be read out from table \ref{TAB:RESULTS}, where again $\SymbolAQPL$ consistently performs best with respect to all measures. The particular differences in performance between the methods vary strongly between the data sets, chosen samplings and contrasts but on average $\SymbolAQPL$ improves on total variation by about 6dB in PSNR and 8\% in SSIM for both contrasts, \cf~table \ref{TAB:RESULTS}.

\subsection{Discussion}
The largest improvements were obtained for the $T_1$-weighted \NameSheppLogan\linebreak and both contrasts from \NameBrainWeb. We attribute this to the higher level of detail in $T_1$ than in $T_2$ version of \NameSheppLogan, which in turn results in $T_1$ having higher total variation than $T_2$. While the quantitative results for \NameNinon~do not indicate much improvement, visual inspection corroborates the increase in image quality. This might be due to the fact that the noisy reconstructed MRI images have been taken as ground truth and in this case the similarity measure do not match human perception.
\section{Conclusions}
In this paper we extended total variation to accommodate the structural a priori information available from another contrast in MRI. The structural information can be either purely on the location or on the location and direction of edges. In both cases, the prior is convex so that we can use efficient methods from convex optimization to solve the problem. The numerical results with numerous test cases show that exploiting structural information is beneficial in the reconstruction of highly undersampled MRI. Moreover, utilizing directional information yields not only better defined images but also better reconstruction of smooth structures and fine details.

In the future, we will extend the proposed framework to more than a pair of contrasts and so to exploit the structural similarity of a whole sequence of MRI images. Moreover, we intend to extend our method to joint reconstruction of multiple contrasts.

\section*{Acknowledgements}
The authors would like to thank Felix Lucka and Ivana Drobnjak for helpful discussions. Moreover, we highly appreciate the help of Ninon Burgos and Jonathan Schott who provided the real MRI images.
% %BibTex
% % \section*{References}
% \bibliographystyle{siam}
% \bibliography{input/library}%
% 
%BibLatex
% \defbibheading{myheading}[References]{%
%   \section*{#1}}
% \printbibliography[heading=myheading]
\printbibliography
\end{document}